\documentclass[twoside,11pt]{article}
\usepackage{graphicx}
\usepackage{color}
\usepackage{amsmath,amssymb,amsthm,amsfonts}
\usepackage{graphicx}
\usepackage{color}
\usepackage{multicol}
\usepackage{enumerate}
\usepackage[numbers,sort&compress]{natbib}
\usepackage{mathrsfs}

\theoremstyle{plain}
\newtheorem{theorem}{Theorem}[section]
\newtheorem{lemma}{Lemma}[section]
\newtheorem{proposition}{Proposition}[section]

\theoremstyle{definition}

\newtheorem{remark}{Remark}[section]

\def\R{\mathbb{R}}

\numberwithin{equation}{section}

\begin{document}
\allowdisplaybreaks
\title{\large\bf TRAVELLING WAVE SOLUTIONS OF THE DENSITY-SUPPRESSED MOTILITY MODEL}
%\author{
%{\rm Jing Li\footnote{ E-mail address: matlj@163.com}}\\
%{\it\small  College of Science, Minzu University of China, Beijing, 100081, P.R. China}\\
%\vspace{0.5cm}
%{\rm Rachidi B. Salako\footnote{ E-mail address: rbs0016@anburn.edu}}\\
%{\it\small  Department of Mathematics and Statistics, Auburn University, AL 36849, USA }\\
%\vspace{0.5cm}
%{\rm Zhian Wang\footnote{ Corresponding author. E-mail address: mawza@polyu.edu.hk}}\\
%{\it\small Department of Applied Mathematics, Hong Kong Polytechnic University, Hung Hom, Hong Kong }
%}
%
%\date{}
%
%\maketitle
\author{Jing Li\thanks{College of Science, Minzu University of China, Beijing, 100081, P.R. China
        ({\tt matlj@163.com}).}
\and Zhi-An Wang \thanks{Corresponding author: Department of Applied Mathematics, Hong Kong Polytechnic University, Hung Hom, Kowloon, Hong Kong
        ({\tt mawza@polyu.edu.hk}).}}
\date{}
\maketitle

\begin{abstract}
To understand the ``self-trapping'' mechanism inducing spatio-temporal pattern formations observed in the experiment of \cite{Liu-Science} for bacterial motion, the following density-suppressed motility model
\begin{align*}
\left\{
\begin{array}{ll}
u_t=\Delta(\gamma(v)u)+u(a-bu),\\
v_t=\Delta v+u-v,
\end{array}
\right.
\end{align*}
was proposed in \cite{Liu-Science,Fu-PRL}, where $u(x,t)$ and $v(x,t)$ represent the densities of bacteria and the chemical emitted by the bacteria, respectively; $\gamma(v)$ is called the motility function satisfying $\gamma'(v)<0$ and $a,b>0$ are positive constants accounting for the growth and death rates of bacterial cells.  The analysis of the above system is highly non-trivial due to the cross-diffusion and possible degeneracy resulting from the nonlinear motility function $\gamma(v)$ and mathematical progresses on the global well-posedness and asymptotics of solutions were just made recently. Among other things,  the purpose of this paper is to consider a specialized motility function $\gamma(v)=\frac{1}{(1+v)^m} (m>0)$ and investigate the travelling wave solutions which are genuine patterns observed in the experiment of \cite{Liu-Science}. By ingeniously introducing an auxiliary parabolic problem to which the comparison principle applies and constructing relaxed super- and sub-solutions with spatially inhomogeneous decay rates, we show that there is a number $c^*(a,b,m)>2\sqrt{a}$ with $c^*(a,b,m)\to \infty$ as $m\to 0$ such that the above density-suppressed motility model admits travelling wave solutions $(u,v)(x,t)=:(U,V)(x\cdot \xi -ct)$ in $\R^N$ along the direction $\xi \in S^{N-1}$ for all wave speed $c\in[2\sqrt{a},c^*(a,b,m))$ connecting the equilibrium $(a/b,a/b)$ to $(0,0)$, while positive travelling wave solutions will not exist if $c<2\sqrt{a}$.   As $m \to 0$, $\gamma(v) \to 1$ and our results are well consistent with the relevant results for the well-known Fisher-KPP equation (i.e. the first equation of the above system with $\gamma(v)=1$). We further discuss the selection of wave patterns and wave speeds for given initial value and use numerical simulations to illustrate that both monotone and non-monotone traveling wavefronts exist depending on whether the motility function $\gamma(v)$ changes its convexity at $v=a/b$. Two-dimensional simulations demonstrate that the system can generate outward expanding ring (strip) pattern as observed in the experiment.
\end{abstract}

\vskip5mm

{\small


\ \ \ \ \ \ \  \ \ \ \ \ \ \ \ auxiliary problem, spatially inhomogeneous decay rate
}
\vskip5mm

\section{Introduction}
The reaction-diffusion models can reproduce a wide variety of exquisite spatio-temporal patterns arising in embryogenesis, development and population dynamics due to the diffusion-driven (Turing) instability \cite{Kondo-Miura-Science-2010, Murray-book-2001}. Many of them invoke nonlinear diffusion enhanced by the local environment condition to accounting for population pressure (cf. \cite{mendez2012density}), volume exclusion (cf. \cite{Painter-Hillen-CAMQ, Wang-Hillen}) or avoidance of danger (cf. \cite{Murray-book-2001}) and so on. However the opposite situation where the species will slow down its random diffusion rate when encountering external signals such as the predator in pursuit of the prey \cite{JW-EJAP-2019, Karevia} and the bacterial in searching food \cite{KS1, KS2} has not been considered.  Recently a so-called ``self-trapping'' mechanism was introduced in \cite{Liu-Science} by a synthetic biology approach onto programmed bacterial {\it Eeshcrichia coli} cells which excrete signalling molecules acyl-homoserine lactone (AHL) such that at low AHL level, the bacteria undergo run-and-tumble random motion and are motile, while at high AHL levels, the bacteria tumble incessantly and become immotile due to the vanishing macroscopic motility. Remarkably {\it Eeshcrichia coli} cells formed the outward expanding ring (strip) patterns in the petri dish. To understand the underlying patterning mechanism, the following two-component ``density-suppressed motility'' reaction-diffusion system has been proposed in \cite{Fu-PRL}
\begin{align}\label{1.1}
\left\{
\begin{array}{ll}
u_t=\Delta(\gamma(v)u)+u(a-bu),
\\
\tau v_t=\Delta v+u-v,
\end{array}
\right.
\end{align}
where $u(x,t),v(x,t)$ denote the bacterial cell density, concentration of acyl-homoserine lactone (AHL) at position $x$ and time $t$, respectively. The first equation of \eqref{1.1} describes the random motion of bacterial cells with an AHL-dependent motility coefficient $\gamma(v)$ and logistic cell growth with growth rate $a>0$ and death rate $b>0$. The second equation of \eqref{1.1} describes the  diffusion, production and turnover of AHL with $\tau\in \{0,1\}$.
%In \cite{Fu-PRL}, $\tau=1$ and for the convenience of statement in the sequel we also allow $\tau=0$.
The striking feature of the system \eqref{1.1} is that the bacterial diffusion rate is a function $\gamma(v)$ depending on an external signal density $v$, which satisfies $\gamma'(v)<0$ accounting for the repressive effect of AHL concentration on the bacterial motility (cf. \cite{Liu-Science}). This monotone decreasing property of $\gamma(v)$ distinguishes the nonlinear diffusion in \eqref{1.1} from other cross-diffusion systems (cf. \cite{Lou}) where the diffusion of species is increasing with respect to density due to population pressure. We remark that the system \eqref{1.1} was originally given in the supplementary material of \cite{Liu-Science} and formally analyzed in \cite{Fu-PRL}.

\begin{figure}[!htbp]
\centering
\includegraphics[width=14cm]{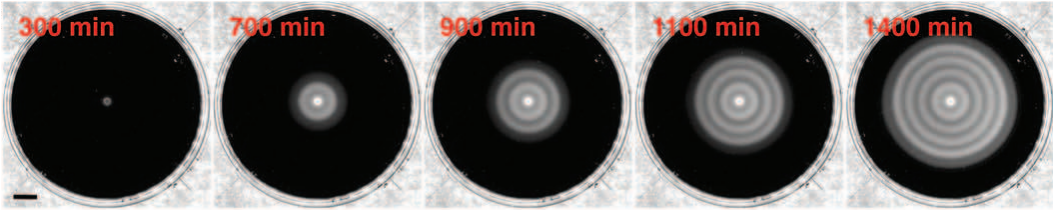}

\caption{Time-lapsed photographs of spatiotemporal patterns formed by the engineered {\it Eeshcrichia coli} strain CL3 (see details in \cite{Liu-Science}). The figure is taken from the Figure 1 in \cite{Liu-Science} for illustration.}
\label{fig0}
\end{figure}

Expansion of the Laplacian term $\Delta (\gamma(v)u)=\nabla\cdot(\gamma(v)\nabla u+u\gamma'(v)\nabla v)$ in the first equation of \eqref{1.1} indicates that the motility function $\gamma(v)$ generates a cross-diffusion effect, and the decay property $\gamma'(v)<0$ may lead to degenerate diffusion making the analysis highly nontrivial. Therefore not many mathematical result have been available to \eqref{1.1} which has received attentions in recent years. When the system \eqref{1.1} is considered in a bounded domain $\Omega$ with Neumann boundary conditions, the following results are obtained in the literature.

\begin{itemize}
\item[(C1)](With cell growth: $a=b>0$) Firstly the global existence and large time behavior of solutions was established in \cite{JKW-SIAP-2018} where it was shown that the system \eqref{1.1} with $\tau=1$ has a unique global classical solution in two dimensions under the following assumptions on the motility function $\gamma(v)$:

\begin{itemize}
\item[(H0)] $\gamma(v)\in C^3([0,\infty)), \gamma(v)>0\ \mathrm{and}~\gamma'(v)<0 $, $\lim\limits_{v \to \infty}\gamma(v)=0$ and $\lim\limits_{v \to \infty}\frac{\gamma'(v)}{\gamma(v)}$ exists.
\end{itemize}

\vspace{0.2cm}
 Moreover, the constant steady state $(1,1)$ of \eqref{1.1} is proved to be globally asymptotically stable if
$
a=b>\frac{K_0}{16}$ where $ K_0=\max\limits_{0\leq v \leq \infty}\frac{|\gamma'(v)|^2}{\gamma(v)}$. Later the global existence result was extended to  higher dimensions ($n\geq 3$)  for large $a>0$ in \cite{WW-JMP-2019}. Recently the similar results have been obtained for \eqref{1.1} with $\tau=0$ in \cite{Fujie-Jiang2, JW-DCDSB-2020} without the condition $\lim\limits_{v \to \infty}\frac{\gamma'(v)}{\gamma(v)}$ in (H0). On the other hand, for small $a>0$, the existence/{\color{black}nonexistence} of nonconstant steady states  of \eqref{1.1} was rigorously established under certain conditions in \cite{MPW-PD-2019} and the periodic pulsating wave was analytically approximated by the multi-scale analysis. When $\gamma(v)$ is a piecewise constant function, the dynamics of discontinuity interface was studied in \cite{SIK-EJAP-2019} and existence of discontinuous traveling wave solutions was established in \cite{Lui}.
\item[(C2)](Without cell growth: $a=b=0$) It turns out the dynamics of \eqref{1.1} with $a=b=0$ are very different from the case $a=b>0$ (with cell growth). With a specialized motility function $\gamma(v)=c_0/v^{k}(k>0)$, the global existence of classical solutions of \eqref{1.1} with $\tau=1$ in any dimensions was established in \cite{YK-AAM-2017} for small $c_0>0$. This smallness assumption on $c_0$ was removed later for the parabolic-elliptic case (i.e. \eqref{1.1} with $\tau=0$) with $0<k<\frac{2}{n-2}$ in \cite{AY-Nonlinearity-2019}. If $\gamma(v)$ decays algebraically and $1\leq n\leq 3$, the global existence of weak solutions of \eqref{1.1} with $\tau=1$ with large initial data was established in \cite{Kim-NARWA-2019}. However the solution of \eqref{1.1} may blow up if $\gamma(v)$ has a faster decay rate. For example, if $\gamma(v)=e^{-\chi v}$, by constructing a Lyapunov functional, it was proved in \cite{JW-PAMS-2019} that  there exists a critical mass $m_*=\frac{4\pi}{\chi}$ such that the solution of \eqref{1.1} with $\tau=1$ exists globally with uniform-in-time bound if $\int_\Omega u_0dx<m_*$ while blows up if $\int_\Omega u_0dx>m_*$ in two dimensions, where $u_0$ is the initial value of $u$. The result of \cite{JW-PAMS-2019} was further refined in \cite{Fujie-Jiang1} by showing that the blow-up time is infinite. When $\gamma(v)$ has both positive lower and upper bounds, the global existence of classical solutions in two dimensions was proved in \cite{TW-M3AS-2017}. Very recently the existence/nonexistence of non-constant stationary solutions as well as pattern formation were explored in \cite{Xu-Wang-2020} via the global bifurcation theory and weak-strong solutions of \eqref{1.1} with $\tau=1$ in any dimensions was explored in \cite{Burger}.
\end{itemize}

As recalled above, the existing results for \eqref{1.1} are confined to the global well-posedness, asymptotic behaviors of solutions and stationary solutions (pattern formation). However the traveling wave solutions, which are  genuinely relevant to the experiment observation of \cite{Liu-Science}, are not investigated mathematically except for a special case that $\gamma(v)$ is piecewise constant. When $\gamma(v)$ is a constant, equations of \eqref{1.1} are decoupled each other and the first equation becomes the well-known Fisher-KPP equation - a benchmark model for the study of traveling wave solutions of reaction-diffusion equations \cite{Murray-book-2001}. However, once $\gamma(v)$ is non-constant, \eqref{1.1} becomes a coupled system with cross-diffusion and the study of traveling wave solutions drastically becomes difficult. The purpose of this paper is to make some progress towards this direction and explore the existence of traveling wave solutions to \eqref{1.1} with allowable wave speeds. To be specific and simple, we consider the following motility function
\begin{equation}\label{motility}
\gamma(v)=\frac{1}{(1+v)^m}, \ m>0
\end{equation}
which fulfills the condition (H0). However our argument can be directly extend to other forms of motility function satisfying (H0), for instance $\gamma(v)=e^{-\chi v}$. But the calculations and conditions ensuring the existence of traveling wave solutions may be different.

To put things in perspective, we rewrite \eqref{1.1} as
\begin{align}\label{1.1n}
\left\{
\begin{array}{ll}
u_t=\nabla\cdot(\gamma(v)\nabla u+u\gamma'(v)\nabla v)+u(a-bu),
\\
v_t=\Delta v+u-v,
\end{array}
\right.
\end{align}
which is a Keller-Segel type chemotaxis model proposed in \cite{KS1} with growth.  For the classical chemotaxis-growth system
\begin{align}\label{1.20}
\left\{
\begin{array}{ll}
u_t=\nabla \cdot (\nabla u-\chi u\nabla v)+u(a-bu),
\\
\tau v_t=\Delta v+u-v,
\end{array}
\right.
\end{align}
travelling wave solutions are investigated in a series of works \cite{SS,SS1,SS2} for both cases $\tau=0$ and $\tau=1$, where $\chi>0$ denotes the chemotactic coefficient. The existence of traveling wave solutions with minimal wave speed depending on $a$ and $\chi$ was obtained, the asymptotic wave speed as $\chi\to0$ as well as the spreading speed were examined in details in \cite{SS,SS1,SS2,SS3} where the major tool used therein to prove the existence of traveling wave solutions  is the parabolic comparison principle. Except traveling wave solutions, the chemotaxis-growth system \eqref{1.20} can also drive other complex patterning dynamics (cf. \cite{KWA, MOW,PH-PhyD}). When the volume filling effect in considered in \eqref{1.20} (i.e. $\chi u\nabla v$ is changed to $\chi u(1-u)\nabla v$), the traveling wave solutions with minimal wave speed were shown to exist in \cite{ou-yuan} for small chemotactic coefficient $\chi>0$. For the original singular Keller-Segel system generating traveling waves without cell growth, we refer to \cite{KS3, Li14, wang12} and references therein.
In contrast to the classical chemotaxis-growth system \eqref{1.20}, both diffusive and chemotactic coefficients in the system \eqref{1.1n} are non-constant. This not only makes the analysis more complex, but also make the parabolic comparison principle inapplicable due to the nonlinear diffusion. In this paper, we shall develop some new ideas (see details in section 2) to tackle the various difficulties induced by the nonlinear motility function $\gamma(v)$ and establish the existence of traveling wave solutions to \eqref{1.1}.

%by using a space-dependent eigenvalue to apply the comparison  principle to get solutions of an axillary parabolic equation on a specified space with exact asymptotic decay rate, based on which the existence of traveling wave solutions to \eqref{1.1} can be obtained as a fixed point of the large-time limits of the axillary parabolic equation. To be precise, we shall consider

The rest of this paper is organized as follows. In section 2, we state our main results on the existence/non-existence of traveling wave solutions to \eqref{1.1} for $(x, t)\in \mathbb{R}^N \times [0,\infty)$ and sketch the proof strategies. In section 3, we derive some preliminary results that will be used in the subsequent sections. In section 4, we construct and study some auxiliary problems connecting to our problem. In section 5, we prove our main theorems via Schauder fixed point theorem and compactness argument based on the results in previous sections. In final section 6, we discuss the possible selection of wave profiles/speeds and use numerical simulations to illustrate the traveling wave patterns.

\section{Main results and proof strategies} We shall establish the existence of traveling wave solutions and wave speed, and explore  how the density-suppressed motility influences  travelling wave profiles and ``the minimal wave speed". In the spatially homogeneous situation the steady states are $(0,0)$ and $(a/b,a/b)$, which are respectively unstable (saddle point) and stable node. This suggests that we should look for travelling
wavefront solutions to \eqref{1.1} connecting $(a/b,a/b)$ to $(0,0)$. Moveover negative $u$ and $v$ have no physical meanings with what we have in mind in the sequel.

A nonnegative solution $(u(x,t),v(x,t))$ is called a travelling wave solution of \eqref{1.1} connecting $(a/b,a/b)$ to $(0,0)$ and propagating in the direction $\xi\in S^{N-1}$ with speed $c$ if it is of the form
$$(u(x,t),v(x,t))=(U(x\cdot\xi-ct),V(x\cdot\xi-ct))=: (U(z), V(z))$$
satisfying the following equations
\begin{align}\label{1.2}
\left\{
\begin{array}{ll}
(\gamma(V)U)''+cU'+U(a-bU)=0,\\
V''+cV'+U-V=0
\end{array}
\right.
\end{align}
and
\begin{equation}\label{bc}
(U(-\infty),V(-\infty))=(a/b,a/b), \ (U(+\infty),V(+\infty))=(0,0).
\end{equation}
where $'=\frac{d}{dz}$. In this paper, we proceed to find the constraints on the parameters to exclude the spatial-temporal pattern formation and guarantee the existence of travelling wave solutions connecting the two constant steady states.

Denote
\begin{align}
b^*(m,a)=3m\left(1+2\sqrt{\frac{a}{1+a}}\left(2+\sqrt{1+\frac1m}\right)\right)\label{b^*}
\end{align}
and
\begin{align}
c^*(a,b,m)=\frac{(b-3m)\sqrt{1+a}}{3m}-2\sqrt{a}\left(1+\sqrt{1+\frac1m}\right).\label{c^*}
\end{align}
For any $b\geq b^*(m,a)$, it can be easily verified that $c^*(a,b,m)\geq2\sqrt{a}$ and the set $[2\sqrt{a},c^*(a,b,m)]$ is non-empty. We obtain the following theorem.
\begin{theorem}\label{main-thm1}
Let $\gamma(v)$ be given in \eqref{motility}. Then for any $b\geq b^*(m,a)$ and $c\in[2\sqrt{a},c^*(a,b,m)]$, the system \eqref{1.1} has a travelling wave solution $(u(x,t),v(x,t))=(U(x\cdot\xi-ct),V(x\cdot\xi-ct))$ with speed $c$ in the direction $\xi \in S^{N-1}$ for all $(x,t)\in \R^N\times[0,+\infty)$, satisfying
\begin{equation}\label{am}
\lim_{z\to+\infty}\frac{U(z)}{e^{-\lambda z}}=1,\quad\lim_{z\to+\infty}\frac{V(z)}{e^{-\lambda z}}=\frac{1}{1+a}
\end{equation}
with $\lambda=\frac{c-\sqrt{c^2-4a}}{2}$ and
\begin{align*}
\liminf_{z\to-\infty}U(z)>0\quad\hbox{and}\quad\liminf_{z\to-\infty}V(z)>0.
\end{align*}
Moreover, if
\begin{align}
\mathcal{K}(m,a)=m\sqrt{\frac{a(1+a)}{m(m+1)}}\left(\sqrt{\frac{a(1+a)}{m(m+1)}}+1\right)^m<1,\label{1.12}
\end{align}
we have
$$
\lim_{z\to-\infty}U(z)=\lim_{z\to-\infty}V(z)=a/b
$$
and
$$\lim_{z\to\pm\infty}U'(z)=\lim_{z\to\pm\infty}V'(z)=0.$$
\end{theorem}

\begin{remark}
Theorem 1.1 and Theorem 1.2 imply that $c=2\sqrt{a}$ is the minimal wave speed same as the one for the classical Fisher-KPP equation and irrelevant to the decay rate of the motility function. Different from the Fisher-KPP equation, the wave speed has an upper bound $c^*(a,b,m)$ induced by  the density-suppressed motility. As $m\to0$, $\gamma(v) \to 1$ and the equation for $u$ becomes the classical Fisher-KPP equation. From the definition of $b^*(m,a)$ and $c^*(a,b,m)$, it holds that
$$\lim_{m\to0} b^*(m,a)=0\ \ \mathrm{and} \ \ \lim_{m\to0}c^*(a,b,m)=+\infty$$
which agree with the results for the classical Fisher-KPP equation.

Moreover it is straightforward to check that $\lim\limits_{m\to0}\mathcal{K}(m,a)=0$ and $\lim\limits_{a\to0}\mathcal{K}(m,a)=0,$ which implies the condition \eqref{1.12} can be ensured when $m>0$ or $a>0$ is small.
\end{remark}
\begin{theorem}\label{main-thm2}
For $c<2\sqrt{a}$, there is no travelling wave solution $(u(x,t),v(x,t))=(U(x\cdot\xi-ct),V(x\cdot\xi-ct))$ of \eqref{1.1} connecting the constant solutions $(a/b,a/b)$ and $(0,0)$ with speed $c$.
\end{theorem}

%\begin{remark}
%Theorem 1.1 and Theorem 1.2 show that $c=2\sqrt{a}$ is the ``minimal wave speed", which means that the density-suppressed motility does not change the minimal wave speed for classical Fisher-KPP equation.
%\end{remark}

\noindent {\bf Proof strategies}. Since the model \eqref{1.1} is a cross diffusion system, see also \eqref{1.1n}, many classical tools proving the existence of traveling waves such as phase plane analysis, topological methods and bifurcation analysis (cf. \cite{Volpert}), among others, become infeasible. Motivated from excellent works of Salako and Shen \cite{SS,SS1} for the chemotaxis-growth model \eqref{1.20} by constructing super- and sub-solutions and proving the existence of traveling wave solutions as the large-time limit of solutions in the moving-coordinate system based on the parabolic comparison principle, we plan to achieve our goals in a similar spirit. However substantial differences exist between the models \eqref{1.1} and \eqref{1.20}. The nonlinear motility function $\gamma(v)$ in \eqref{1.1} refrains us employing the parabolic comparison principle and constructing super- and sub-solutions with the same decay rate at the far field, which are crucial ingredients used for \eqref{1.20}  in \cite{SS,SS1}. In this paper, we develop two innovative ideas to overcome these barriers. First we introduce an auxiliary parabolic problem \eqref{A2} with constant diffusion to which the method of super- and sub-solutions applies (see section 4.1). This auxiliary problem subtly bypasses the barriers induced by the nonlinear diffusion but its time-asymptotic limit yields a solution to an elliptic problem \eqref{1.4} whose fixed points indeed correspond to solutions to \eqref{1.2} - namely traveling wave solutions to our concerned system \eqref{1.1} (see section 4.2). Second we construct a sequence of relaxed sub-solution $\underline{U}_n(x)$ for any $n>1$ with a spatially inhomogeneous decay rate $\theta_1(x)$ which approaches to the constant decay rate of the super-solution $\overline{U}(x)$ as $x \to +\infty$ (see section 3.1). With them we use the method of super- and sub-solutions to construct solutions to the auxiliary parabolic problem \eqref{A2} in appropriate function space and manage to show its time-asymptotic limit problem has a fixed point. This is a fresh idea substantially different from the works \cite{SS,SS1} where the super- and sub-solutions were directly constructed with the same decay rates by taking the advantage of constant diffusion.

We divide the proof of Theorem \ref{main-thm1} into four steps. In step 1, we construct an auxiliary parabolic problem \eqref{A2} with constant diffusion and prove its global boundedness uniformly in time (see Proposition \ref{Pro2.2}) by the method of super- and sub-solutions. In step 2, we show that the limit of global solutions to \eqref{A2} as $t \to \infty$ yields a semi-wavefront solution to an elliptic problem \eqref{1.4} with some compactness argument (see Proposition \ref{proposition2.1}). In step 3, we show that the solution obtained in Step 2 satisfies the boundary condition \eqref{bc} by direct estimates under some constraints on $m$ and $a$ (see Proposition \ref{Pro2.3}), which hence warrants that the semi-wavefront solution is indeed a wavefront solution in $\R$. Finally in step 4, we use the Schauder's fixed point theorem to prove that \eqref{1.4} has a fixed point which gives a solution to \eqref{1.2} in $\R$ satisfying \eqref{bc} (see section 5.1), where the trick of utilizing relaxed sub-solution $\underline{U}_n(x)$ with spatially inhomogeneous decay rate is critically used to obtain the continuity of solution map. Theorem \ref{main-thm2} is proved directly by an argument of contradiction.

\section{Preliminary results}
In this section we introduce some notations/definitions and list some basic facts which will be used in our subsequent analysis. In particular, the construction of relaxed super and sub-solutions with spatially inhomogeneous decay rates will be presented in this section as a preparation for the analysis in section 4.

\subsection{Super and sub-solutions with spatially inhomogeneous decay rates}
\begin{figure}[htbp]
\centering
\includegraphics[width=7cm,angle=0]{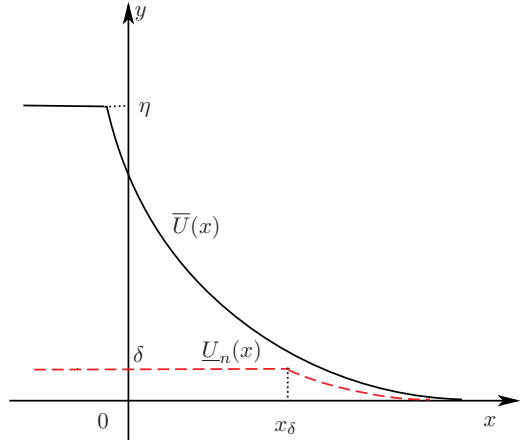}
\caption{A schematic of functions $\overline{U}(x)$ and $\underline{U}_n(x)$, where the solid black line represents $\overline{U}(x)$ and the dashed red line represents $\underline{U}_n(x)$.}
\label{fig}
\end{figure}

For $c\geq 2\sqrt{a}$, define
\begin{align}
\lambda:=\frac{c-\sqrt{c^2-4a}}{2}\quad\hbox{and}\quad\theta_1(x):=\frac{c-\sqrt{c^2-4a\Big(1+\frac{e^{-\lambda x}}{1+a}\Big)^{-m}}}{2\Big(1+\frac{e^{-\lambda x}}{1+a}\Big)^{-m}}\quad\forall x\in\mathbb R,\label{1.21}
\end{align}
for which
\begin{align}
\lambda^2-c\lambda+a=0,\qquad\Big(1+\frac{e^{-\lambda x}}{1+a}\Big)^{-m}\theta_1^2(x)-c\theta_1(x)+a=0\quad\forall x\in\mathbb R\label{1.23}
\end{align}
and
\begin{align}
\lim_{x\to+\infty}\theta_1(x)=\lambda,\quad 0<\theta_1(x)<\lambda\leq\sqrt{a}\quad\forall x\in\mathbb R.\label{1.22}
\end{align}
Choose
\begin{align}\label{1.29}
\theta_2(x):=\left\{\begin{array}{ll}
\theta_1(x)+\lambda/4,  \quad c=2\sqrt{a},
\\
\theta_1(x)+\lambda/k_0,\quad c>2\sqrt{a}
\end{array}\right.\quad\forall x\in\mathbb R
\end{align}
with $k_0>\max\left\{\frac{2\lambda}{c-2\lambda},2\right\}$. Then
\begin{align}
\theta_2(x)\in\left(\theta_1(x),\theta_1(x)+\frac\lambda 2\right)\quad\forall x\in\mathbb R\label{1.24}
\end{align}
and there exists $x_0\in \mathbb{R}$ such that
\begin{align}
\theta_2(x)<2\theta_1(x)\quad\hbox{for}\; x>x_0.\label{1.30}
\end{align}
Define two functions:
\begin{align}
\overline{U}(x):=\min\{e^{-\lambda x},\eta \}\quad\forall x\in\mathbb R\label{1.25}
\end{align}
and
\begin{align}\label{1.26}
\underline{U}_n(x):=\left\{
\begin{array}{ll}
\delta, & \quad x\leq x_\delta,
\\
d_ne^{-\theta_1(x)x}+d_0e^{-\theta_2(x) x},& \quad x>x_\delta,
\end{array}
\right.
\end{align}
where $\delta$ is chosen sufficiently small,
$x_\delta>0$ is the unique positive solution of the equation $d_ne^{-\theta_1(x)x}+d_0e^{-\theta_2(x) x}=\delta$,
\begin{align}\label{dn}
d_n:=1-\frac1n\quad\hbox{with}\; 2\leq n\in\mathbb N,\qquad
d_0:=\left\{\begin{array}{ll}
1, &c=2\sqrt{a},
\\
-1, &c>2\sqrt{a}
\end{array}\right.
\end{align}
and
\begin{align}\label{2.8}
\eta :=\frac{b-m\left(1+\frac{c}{\sqrt{1+a}}\right)-\sqrt{\left(b-m\left(1+\frac{c}{\sqrt{1+a}}\right)\right)^2-\frac{4m(m+1)a}{1+a}}}{\frac{2m(m+1)}{1+a}}.
\end{align}
Noticing that $d_n\in(0,1)$ and $$\lim_{x\to+\infty}e^{(\theta_1(x)-\lambda)x}=1,$$ which will be verified in Lemma \ref{Lemma2.1},
we can choose sufficiently small $\delta$, with which $x_\delta$ is large enough such that for all $x\in\mathbb R$, $$0<\underline{U}_n<\overline{U}\leq \eta .$$
We note that the functions $\overline{U}(x)$ and $\underline{U}_n(x)$ will be essentially used later as the super- and sub-solutions of an auxiliary problem we introduce in section 4. A schematic of $\overline{U}(x)$ and $\underline{U}_n(x)$ is plotted in Fig.\ref{fig}. Note that the coefficient $d_n$ ($n\geq 2$) and $d_0$ determines the amplitude of $\underline{U}_n(x)$ and $\theta_1(x)$ determine the decay of $\underline{U}_n(x)$ for large $x>x_\delta$. This is a new ingredient developed in this paper to settle the difficulty of analysis caused by the nonlinear motility function $\gamma(v)$.

Denote
$$C_{\mathrm{unif}}^b(\mathbb R):=\{u\in C(\mathbb R)|\ u\;\hbox{is uniformly continuous in}\; \mathbb R \ \mathrm{and} \ \sup_{z\in\mathbb R}|u(z)|<+\infty\},$$
which is equipped with the norm $$\|u\|=\sup_{z\in\mathbb R}|u(z)|.$$
Define the function space
\begin{equation}\label{fs}
\mathscr{E}_n:=\{u\in C_{\mathrm{unif}}^b(\mathbb R)|\underline{U}_n\leq u\leq \overline{U}\}
\end{equation}
and
\begin{align*}
\mathscr{X}_0:=\bigcap_{n>1}\mathscr{E}_n.
\end{align*}
To find solutions of \eqref{1.2} in $\mathscr{X}_0$, we need the following Lemmas.

\begin{lemma}\label{Lemma2.1}
Let $\lambda$ and $\theta_1(x)$ be defined in \eqref{1.21}. Then it follows that
\begin{align}
\lim_{x\to+\infty}e^{(\theta_1(x)-\lambda)x}=1.\label{1.47}
\end{align}
Moreover, for sufficiently small $\delta>0$, if $x>x_\delta$, then for $c=2\sqrt{a}$,
\begin{align}
0<\theta_1'(x)\leq 2K_1e^{-\frac\lambda2 x} \ \mathrm{and}\  -\lambda K_1e^{-\frac\lambda2 x}\leq  \theta_1''(x)<0 \label{1.48}
\end{align}
with $K_1=\frac{a}{2}\sqrt{\frac{m}{1+a}}$; while for $c>2\sqrt{a}$,
\begin{align}
0<\theta_1'(x)\leq 2K_2e^{-\lambda x} \ \mathrm{and} \ -2\lambda K_2e^{-\lambda x}\leq \theta_1''(x)<0\label{1.50}
\end{align}
with
$K_2=\frac{4a^2m\lambda}{\left(c+\sqrt{c^2-4a}\right)^2\sqrt{c^2-4a}(1+a)}.$
\end{lemma}

\proof
In the sequel, for notational simplicity, under $c\geq 2\sqrt{a}$, we introduce the following notations
\begin{eqnarray*}
\begin{cases}
\phi(x):=\displaystyle 1+\frac{e^{-\lambda x}}{1+a},\\[2mm]
\rho(\phi(x)):=\sqrt{c^2-4a\phi^{-m}(x)},\\[2mm]
%\rho(\phi(x)):=\rho_1(x;c,a,m)=\left(c+\sqrt{c^2-\frac{4a}{\phi(x)^m}}\right)^2,\\[2mm]
h(\phi(x)):=\displaystyle \frac{2a}{c+\rho(\phi(x))}=\displaystyle \frac{2a}{c+\sqrt{c^2-4a\phi^{-m}(x)}}.
\end{cases}
\end{eqnarray*}
Then from the definition of $\theta_1(x)$, we  have
$$\theta_1(x)=\frac{c-\sqrt{c^2-4a\phi^{-m}(x)}}{2\phi^{-m}(x)}=\frac{2a}{c+\sqrt{c^2-4a\phi^{-m}(x)}}=h(\phi(x)).$$
%Denote $h(s):=\frac{2a}{c+\sqrt{c^2-\frac{4a}{\phi(x)^m}}}$ and $s(x):=\frac{e^{-\lambda x}}{1+a}$,
With simple calculation, we find
\begin{align}\label{2.14}
h'(\phi(x))=\frac{-4a^2m}{\rho(\phi(x))(c+\rho(\phi(x)))^2\phi^{m+1}(x)},\quad\phi'(x)=\frac{-\lambda e^{-\lambda x}}{1+a}
\end{align}
and then
\begin{align}
\theta_1'(x)&=(h(\phi(x)))'=h'(\phi)\phi'(x)=\frac{4a^2m\lambda e^{-\lambda x}}{\rho(\phi(x))(c+\rho(\phi(x)))^2\phi^{m+1}(x)(1+a)}>0.\label{1.55}
\end{align}
When $c=2\sqrt{a}$, it has that $\lim\limits_{x\to+\infty}\phi(x)=1$ and $\lim\limits_{x\to+\infty}\rho(\phi(x))=0$. By L'Hopital's rule, we have
\begin{align}
\lim_{x\to+\infty}\bigg|\frac{e^{-\frac{\lambda}{2} x}}{\rho(\phi(x))}\bigg|^2&=\lim_{x\to+\infty}\frac{e^{-\lambda x}}{c^2-4a\phi^{-m}(x)}\nonumber
\\
&=\lim_{x\to+\infty}\frac{e^{-\lambda x}}{4a\phi^m(x)-4a}\lim_{x\to+\infty}\phi^m(x)\label{1.57}
\\
&=\lim_{x\to+\infty}\frac{1+a}{4ma\phi^{m-1}(x)}=\frac{1+a}{4ma}.\nonumber
\end{align}
Then it can be easily verified that
$$\lim_{x\to+\infty}\frac{4a^2m\lambda e^{-\frac{\lambda}{2}x}}{\rho(\phi(x))(c+\rho(\phi(x)))^2\phi^{m+1}(x)(1+a)}=\frac{a}{2}\sqrt{\frac{m}{1+a}}:=K_1,$$
from which and \eqref{1.55}, by choosing sufficiently small $\delta>0$, we can find $x_\delta>0$ such that for $x>x_\delta$, there holds that
 $$0<\theta_1'(x)\leq 2K_1e^{-\frac{\lambda}{2}x}.$$
When $c>2\sqrt{a}$, it has that $\lim\limits_{x\to+\infty}\phi(x)=1$ and $\lim\limits_{x\to+\infty}\rho(\phi(x))=\sqrt{c^2-4a}$. It can be directly checked that
\begin{align*}
\lim_{x\to+\infty}\frac{4a^2m\lambda }{\rho(\phi(x))(c+\rho(\phi(x)))^2\phi^{m+1}(x)(1+a)}=\frac{4a^2m\lambda}{\sqrt{c^2-4a}(c+\sqrt{c^2-4a})^2(1+a)}:=K_2.
\end{align*}
Then from \eqref{1.55}, by choosing sufficiently small $\delta>0$, for $x>x_\delta$,
we obtain $$0<\theta_1'(x)\leq  2K_2e^{-\lambda x}.$$ The first parts of \eqref{1.48} and \eqref{1.50} are proved.

On the other hand,
by L'H\^{o}pital's rule, using \eqref{2.14} and \eqref{1.57}, for $c=2\sqrt{a}$, we obtain
\begin{align*}
&\lim_{x\to+\infty}\left(h(\phi(x))-h(1)\right)x
\\
=&\lim_{x\to+\infty}h'(\phi(x))\frac{\lambda x^2e^{-\lambda x}}{1+a}
\\
=&\lim_{x\to+\infty}\frac{-4a^2m\lambda x^2 e^{-\lambda x}}{\rho(\phi(x))(c+\rho(\phi(x)))^2\phi^{m+1}(x)(1+a)}
\\
=&-\frac{4a^2m\lambda}{c^2(1+a)}\lim_{x\to+\infty}\frac{e^{-\frac\lambda2 x}}{\rho(\phi(x))}\lim_{x\to+\infty}\frac{x^2}{e^{\frac\lambda2 x}}
\\
=&-\frac{am\lambda}{1+a}\sqrt{\frac{1+a}{4ma}}\lim_{x\to+\infty}\frac{x^2}{e^{\frac\lambda2 x}}=0.
\end{align*}
While for $c>2\sqrt{a}$, we obtain
\begin{align*}
\lim_{x\to+\infty}\left(h(\phi(x))-h(1)\right)x=&\lim_{x\to+\infty}h'(\phi(x))\frac{\lambda x^2e^{-\lambda x}}{1+a}
=\frac{\lambda h'(1)}{1+a}\lim_{x\to+\infty}\frac{x^2}{e^{\lambda x}}=0.
\end{align*}
Summing up, for $c\geq2\sqrt{a}$, we obtain
$$\lim_{x\to+\infty}e^{(\theta_1(x)-\lambda)x}=e^{\lim\limits_{x\to+\infty}(h(\phi(x))-h(1))x}=1$$
and then \eqref{1.47} follows.

Now we turn to the estimate of $\theta_1''(x)$. Noticing that
$
h''(\phi)=h'(\phi)(\ln(- h'(\phi)))'\label{2.15}
$
and
\begin{align}
(\ln(-h'(\phi)))'=&\left[\ln(4a^2m)-2\ln(c+\rho(\phi))-\ln \rho(\phi)-(m+1)\ln\phi\right]'\nonumber
\\
=&\frac{-4am}{\rho(\phi)(c+\rho(\phi))\phi^{m+1}}-\frac{2am}{(\rho(\phi))^2\phi^{m+1}}-\frac{m+1}{\phi},\nonumber %\label{2.24}
\end{align}
which together with \eqref{2.14} and the fact that $\phi'(x)=\frac{-\lambda e^{-\lambda x}}{1+a}$ and $\phi''(x)=\frac{\lambda^2 e^{-\lambda x}}{1+a}$ implies
\begin{eqnarray*}\label{1.56}
\begin{aligned}
\theta_1''(x)
&=(h(\phi(x))''=(h'(\phi(x))\phi'(x))'=h'(\phi(x))\phi''(x)+h''(\phi(x))(\phi'(x))^2
\\
&=h'(\phi(x))[\phi''(x)+(\ln (-h'(\phi)))'(\phi'(x))^2]
\\
&=\frac{-4a^2m\lambda^2e^{-\lambda x}}{\rho(\phi(x))(c+\rho(\phi(x)))^2\phi(x)^{m+1}(1+a)}\left\{1-\frac{e^{-\lambda x}}{1+a}\left(\frac{4am}{\rho(\phi(x))(c+\rho(\phi(x)))\phi(x)^{m+1}}\right.\right.
\\
&\quad\left.\left.+\frac{2am}{(\rho(\phi(x)))^2\phi(x)^{m+1}}+\frac{m+1}{\phi(x)}\right)\right\}.
\end{aligned}
\end{eqnarray*}
For $c=2\sqrt{a}$, then $\lambda=\sqrt{a}$, from \eqref{1.57}, it can be verified that
$$\lim_{x\to+\infty}\frac{-4a^2m\lambda^2e^{-\frac\lambda2 x}}{\rho(\phi(x))(c+\rho(\phi(x)))^2\phi(x)^{m+1}(1+a)}=-\lambda K_1$$
and
$$\lim_{x\to+\infty}\frac{e^{-\lambda x}}{1+a}\left(\frac{4am}{\rho(\phi(x))(c+\rho(\phi(x)))\phi(x)^{m+1}}
+\frac{2am}{(\rho(\phi(x)))^2\phi(x)^{m+1}}+\frac{m+1}{\phi(x)}\right)=\frac12.$$
By choosing sufficiently small $\delta>0$,we can find a   $x_\delta>0$ such that
$$0>\theta_1''(x)\geq-\lambda K_1e^{-\frac{\lambda}{2}x}, \ \mathrm{for} \ x>x_\delta.$$
While for $c>2\sqrt{a}$, we can check that
$$\lim_{x\to+\infty}\frac{-4a^2m\lambda^2}{\rho(\phi(x))(c+\rho(\phi(x)))^2\phi(x)^{m+1}(1+a)}=-\lambda K_2$$
and
$$\lim_{x\to+\infty}\frac{e^{-\lambda x}}{1+a}\left(\frac{4am}{\rho(\phi(x))(c+\rho(\phi(x)))\phi(x)^{m+1}}
+\frac{2am}{(\rho(\phi(x)))^2\phi(x)^{m+1}}+\frac{m+1}{\phi(x)}\right)=0.$$
Then by choosing sufficiently small $\delta>0$ so as to generate a  $x_\delta>0$, we have
\begin{align*}
0>\theta_1''(x)\geq-2\lambda K_2e^{-\lambda x}, \ \mathrm{for} \ x>x_\delta.
\end{align*}
Then the last parts of \eqref{1.48} and \eqref{1.50} follows. This completes the proof of Lemma \ref{Lemma2.1}.

$\hfill\square$

\subsection{Some priori estimates}
\begin{lemma}\label{Lemma2.2}
For any $u\in\mathscr{E}_n$, denote $V(\cdot;u)$ the solution of
\begin{align}
V''+cV'+u-V=0.\label{1.18}
\end{align}
Then for $c\geq2\sqrt{a}$, we have
\begin{align}
&0<V(x;u)\leq\min\left\{\frac{e^{-\lambda x}}{1+a},\eta \right\}, \label{2.1}\\
&|V'(x;u)|\leq \min\left\{\frac{\eta }{\sqrt{1+a}},\frac{e^{-\lambda x}}{\sqrt{1+a}}\right\} \label{2.2}
\end{align}
for all $x\in\mathbb R$.
\end{lemma}

\proof Denote
\begin{align}
\lambda_1=\frac{-c-\sqrt{c^2+4}}{2},\quad \lambda_2=\frac{-c+\sqrt{c^2+4}}{2}.\label{1.38}
\end{align}
From \eqref{1.38} and the definition of $\lambda$ in \eqref{1.21}, we obtain
\begin{align}
0<\lambda\leq\sqrt{a},\quad\lambda_1<0,\quad\lambda_2>0,\quad\lambda_1+\lambda<0,\quad\lambda_2+\lambda>0\label{1.45}
\end{align}
and
\begin{align}
\lambda_1\lambda_2=-1,\quad \lambda_1+\lambda_2=-c,\quad \lambda^2-c\lambda-1=-(1+a).\label{1.46}
\end{align}
By the variation of constants, the solution of \eqref{1.18} can be expressed as
\begin{align}
V(x;u)=\frac{1}{\lambda_2-\lambda_1}\left(\int_{-\infty}^xe^{\lambda_1(x-s)}u(s)ds+\int_x^{+\infty}
e^{\lambda_2(x-s)}u(s)ds\right).\label{1.17}
\end{align}
Note that $0\leq u\leq\overline{U}=\min\{\eta ,e^{-\lambda x}\}$ since $u\in\mathscr{E}_n$. Then using \eqref{1.45} and \eqref{1.46}, we obtain from  \eqref{1.17} that
\begin{align*}
0\leq V(x;u)\leq&\frac{1}{\lambda_2-\lambda_1}\left(\int_{-\infty}^xe^{\lambda_1(x-s)}e^{-\lambda s}ds+\int_x^{+\infty}
e^{\lambda_2(x-s)}e^{-\lambda s}ds\right)
\\
=&\frac{1}{\lambda_2-\lambda_1}\left(\frac{e^{-(\lambda_1+\lambda)s}\big|_{-\infty}^{x}}{-(\lambda_1+\lambda)e^{-\lambda_1x}}
+\frac{e^{-(\lambda_2+\lambda)s}\big|_{x}^{+\infty}}{-(\lambda_2+\lambda)e^{-\lambda_2x}}\right)
\\
=&\frac{e^{-\lambda x}}{\lambda_2-\lambda_1}\left(\frac{1}{-(\lambda_1+\lambda)}-\frac{1}{-(\lambda_2+\lambda)}\right)
\\
=&\frac{-e^{-\lambda x}}{\lambda_1\lambda_2+(\lambda_1+\lambda_2)\lambda+\lambda^2}
\\
=&\frac{-e^{-\lambda x}}{\lambda^2-c\lambda-1}=\frac{e^{-\lambda x}}{1+a}
\end{align*}
and
\begin{align*}
0\leq V(x;u)\leq&\frac{1}{\lambda_2-\lambda_1}\left(\int_{-\infty}^xe^{\lambda_1(x-s)}\eta ds+\int_x^{+\infty}
e^{\lambda_2(x-s)}\eta ds\right)
\\
=&\frac{\eta }{\lambda_2-\lambda_1}\left(\frac{e^{-\lambda_1 s}\big|_{-\infty}^x}{-\lambda_1e^{-\lambda_1x}}+\frac{e^{-\lambda_2 s}\big|_x^{+\infty}}{-\lambda_2e^{-\lambda_2x}}\right)
\\
=&\frac{\eta }{\lambda_2-\lambda_1}\left(\frac{1}{-\lambda_1}-\frac{1}{-\lambda_2}\right)
\\
=&\frac{-\eta }{\lambda_1\lambda_2}
=\eta.
\end{align*}
Thus \eqref{2.1} follows. On the other hand, differentiating \eqref{1.17} with respect to $x$, we have
\begin{align}
V'(x;u)%&=\frac{1}{\lambda_2-\lambda_1}\left(\int_{-\infty}^xe^{\lambda_1(x-s)}u(s)ds+\int_x^{+\infty}
%e^{\lambda_2(x-s)}u(s)ds\right)'}]
%\\
=\frac{1}{\lambda_2-\lambda_1}\left(\int_{-\infty}^x\lambda_1e^{\lambda_1(x-s)}u(s)ds+\int_x^{+\infty}\lambda_2
e^{\lambda_2(x-s)}u(s)ds\right).\label{2.13}
\end{align}
For $c\geq 2\sqrt{a}$, using \eqref{1.45}, \eqref{1.46} and the facts that
$$
\lambda_2-\lambda_1=\sqrt{c^2+4}\geq2\sqrt{1+a},\quad
c\lambda=\lambda^2+a<2a,
$$
as well as the fact $0\leq u\leq\min\{\eta ,e^{-\lambda x}\}$, we obtain
\begin{align*}
|V'(x;u)|\leq&\frac{1}{\lambda_2-\lambda_1}\left(\int_{-\infty}^x(-\lambda_1)e^{\lambda_1(x-s)}e^{-\lambda s}ds+\int_x^{+\infty}\lambda_2
e^{\lambda_2(x-s)}e^{-\lambda s}ds\right)
\\
=&\frac{1}{\lambda_2-\lambda_1}\left(\frac{-\lambda_1e^{-(\lambda_1+\lambda)s}\big|_{-\infty}^{x}}{-(\lambda_1+\lambda)e^{-\lambda_1x}}
+\frac{\lambda_2e^{-(\lambda_2+\lambda)s}\big|_{x}^{+\infty}}{-(\lambda_2+\lambda)e^{-\lambda_2x}}\right)
\\
=&\frac{e^{-\lambda x}}{\lambda_2-\lambda_1}\left(\frac{\lambda_1}{\lambda_1+\lambda}+\frac{\lambda_2}{\lambda_2+\lambda}\right)
\\
=&\frac{e^{-\lambda x}}{\lambda_2-\lambda_1}\left(\frac{2\lambda_1\lambda_2+(\lambda_1+\lambda_2)\lambda}{\lambda_1\lambda_2+\lambda^2+(\lambda_1+\lambda_2)\lambda}\right)
\\
=&\frac{-e^{-\lambda x}(2+c\lambda)}{\sqrt{c^2+4}(\lambda^2-c\lambda-1)}
\leq\frac{e^{-\lambda x}}{\sqrt{1+a}}
\end{align*}
and
\begin{align*}
|V'(x;u)|\leq&\frac{1}{\lambda_2-\lambda_1}\left(\int_{-\infty}^x(-\lambda_1)e^{\lambda_1(x-s)}\eta ds+
\int_x^{+\infty}\lambda_2e^{\lambda_2(x-s)}\eta ds\right)
\\
=&\frac{\eta }{\lambda_2-\lambda_1}\left(\frac{e^{-\lambda_1 s}\big|_{-\infty}^x}{e^{-\lambda_1x}}-
\frac{e^{-\lambda_2 s}\big|_x^{+\infty}}{e^{-\lambda_2x}}\right)
\\
=&\frac{2\eta }{\lambda_2-\lambda_1}<\frac{\eta }{\sqrt{1+a}},
\end{align*}
from which \eqref{2.2} follows. The Lemma is thus proved.$\hfill\square$

\begin{lemma}\label{Lemma2.3}
For any $u\in\mathscr{E}_n$, denote $V(x;u)$ the solution of
\begin{align*}
V''+cV'+u-V=0.
\end{align*}
Then for sufficiently small $\delta>0$, if $x>x_\delta$, then
\begin{align}
\gamma(V)\theta_1^2(x)-c\theta_1(x)+a\geq0,\label{2.16}
\end{align}
and
\begin{align}
&\gamma(V)\theta_2^2(x)-c\theta_2(x)+a\geq \frac{a}{64} \quad\mathrm{if}\; c=2\sqrt{a},\label{2.17}
\\
&\gamma(V)\theta_2^2(x)-c\theta_2(x)+a\leq-\frac{\lambda(c-2\lambda)}{4k_0}\quad\mathrm{if}\; c>2\sqrt{a}.\label{2.18}
\end{align}
\end{lemma}

\proof Noticing $V(x)\leq\frac{e^{-\lambda x}}{1+a}$ for $x>x_\delta$, we get \eqref{2.16} from the fact that
\begin{align}
\gamma(V)\theta_1^2(x)-c\theta_1(x)+a\geq\Big(1+\frac{e^{-\lambda x}}{1+a}\Big)^{-m}\theta_1^2(x)-c\theta_1(x)+a=0.\label{1.31}
\end{align}
With
\begin{align*}
\lim_{x\to+\infty}\Big(1+\frac{e^{-\lambda x}}{1+a}\Big)^{-m}=1,\quad\lim_{x\to+\infty}\theta_1(x)=\lambda,
\end{align*}
by choosing sufficiently small $\delta>0$, for all $x>x_\delta$, we have
\begin{align}
\frac{15}{16}\lambda\leq\theta_1(x)<\lambda,\quad \gamma(V)=\frac{1}{(1+V)^m}\geq\Big(1+\frac{e^{-\lambda x}}{1+a}\Big)^{-m}\geq\frac{33}{34}.\label{2.19}
\end{align}
For the case $c=2\sqrt{a}$, for which $\lambda=\frac c2$, noticing $\theta_2(x)=\theta_1(x)+\frac{1}{4}\lambda$, using \eqref{2.19}, we get
\begin{align}
\gamma(V)\theta_2^2(x)-c\theta_2(x)+a&=\gamma(V)\theta_1^2(x)-c\theta_1(x)+a+\gamma(V)\left(\frac1{16}\lambda^2+\frac12\lambda\theta_1(x)\right)
-\frac14c\lambda\nonumber
\\
&\geq\gamma(V)\left(\frac1{16}\lambda^2+\frac12\lambda\theta_1(x)\right)
-\frac14c\lambda\label{1.32}
\\
&\geq\frac{33}{34}\left(\frac1{16}\lambda^2
+\frac{15}{32}\lambda^2\right)-\frac{1}{2}\lambda^2=\frac{a}{64}\nonumber
\end{align}
for all $x>x_\delta$, from which \eqref{2.17} follows.

On the other hand, for the case $c>2\sqrt{a}$, for which $\lambda<\frac c2$, noticing
$$\theta_2(x)=\theta_1(x)+\frac{1}{k_0}\lambda\quad\hbox{with}\;k_0>\max\left\{\frac{2\lambda}{c-2\lambda},2\right\},$$ we obtain
$$\frac{1}{k_0}\lambda^2+2\lambda\theta_1(x)-c\lambda\leq\frac12\lambda(2\lambda-c)<0$$
and then
\begin{align}
&\gamma(V)\theta_2^2(x)-c\theta_2(x)+a\nonumber
\\
\leq&\theta_2^2(x)-c\theta_2(x)+a\label{2.22}
\\
=&\left(\theta_1(x)+\frac\lambda {k_0}\right)^2-c\left(\theta_1(x)+\frac\lambda {k_0}\right)+a\nonumber
\\
=&\theta_1(x)^2-c\theta_1(x)+a+\frac1{k_0}\left(\frac{1}{k_0}\lambda^2+2\lambda\theta_1(x)-c\lambda\right)\nonumber
\\
\leq&\theta_1(x)^2-c\theta_1(x)+a+\frac{1}{2k_0}\lambda(2\lambda-c).\nonumber
\end{align}
moreover owing to the fact $\lim_{x\to+\infty}\theta_1(x)=\lambda$, we have
$$\lim_{x\to+\infty}(\theta_1(x)^2-c\theta_1(x)+a)=\lambda^2-c\lambda+a=0.$$
Then choosing $\delta$ sufficiently small, we obtain that for all $x>x_\delta$.
\begin{align}
\theta_1(x)^2-c\theta_1(x)+a\leq\frac{1}{4k_0}\lambda(c-2\lambda).\label{2.21}
\end{align}
Inserting \eqref{2.21} into \eqref{2.22}, we obtain
\begin{align*}
\gamma(V)\theta_2^2(x)-c\theta_2(x)+a
\leq&\frac{1}{4k_0}\lambda(2\lambda-c)<0.
\end{align*}
Thus \eqref{2.18} follows and Lemma \ref{Lemma2.3} is proved.$\hfill\square$

\section{Auxiliary problems}
In this section, we shall investigate some auxiliary problems which act as bridges to our concerned problem.
\subsection{An auxiliary parabolic problem}

In the sequel, for convenience, we use $\gamma'(v)$ and $\gamma''(v)$ to denote the first and second order derivatives of $\gamma(v)$ with respect to $v$, respectively. This should not be confused with $U', V', U'', V''$ where the prime $'$ means the differentiation with respect to $x$. Given $u\in\mathscr{E}_n$, we first consider the following equation
\begin{equation}\label{A1}
V''+cV'+u-V=0
\end{equation}
which, subject to variation of constants, yields
\begin{align}
V:=V(x;u)=\frac{1}{\lambda_2-\lambda_1}\left(\int_{-\infty}^xe^{\lambda_1(x-s)}u(s)ds+\int_x^{+\infty}
e^{\lambda_2(x-s)}u(s)ds\right).\label{2.23}
\end{align}
Now taking $V$ in \eqref{2.23} as a known function, we define
\begin{align*}
F(U,U'):=\frac{1}{\gamma(V)}\left\{\left(2\gamma'(V)V'+c\right)U'+\left[\gamma''(V)|V'|^2+\gamma'(V)(V-U-cV')+a\right]U-bU^2\right\}.
\end{align*}
By $U(x,t;u,\overline U)$, we denote the solution of the following Cauchy problem
\begin{align}\label{A2}
\left\{\begin{array}{ll}
U_t=U''+F(U,U'), \ x \in \mathbb{R}, t>0
\\[1mm]
U(x,0;u,\overline{U})=\overline{U}(x),\ x \in \mathbb{R}.
\end{array}\right.
\end{align}
From Lemma \ref{Lemma2.2} and the definition of $\gamma(\cdot)$, the boundedness of  $\frac{1}{\gamma(V)}$, $\gamma'(V)$, $\gamma''(V)$, $V$, and $V'$ has been guaranteed. Then the comparison principle is applicable to $\eqref{A2}$. By the semigroup theory,  $U$ can be represented as
\begin{align}
U(x,t;u,\overline U)=&e^{t(\Delta-1)}\overline{U}(x)+\int_0^te^{-(t-s)}e^{(t-s)\Delta}(U+F(U,U'))(x,s)ds.
\end{align}
The local existence of solutions to \eqref{A2} can be obtained by the well-known fixed point theorem (cf. see \cite[Theorem 1.1]{SS2}) along with standard parabolic estimates. We omit the details here for brevity and assume that the solution of \eqref{A2} exists in an maximal interval $[0,T)$ for some $T\in (0, \infty]$ with $U(x,0;u,\overline{U})>0$ for $x\in\mathbb R$. Then the comparison principle for \eqref{A2} implies that $U(x,t;u,\overline{U})>0$ for all $(x,t)\in\mathbb{R}\times [0,T)$.

\begin{proposition}\label{Pro2.2}
For any $b\geq b^*(m,a)$ and $c\in[2\sqrt{a},c^*(a,b,m)]$ with $b^*$, $c^*$ defined in \eqref{b^*} and \eqref{c^*}, there exists $\delta>0$ such that for any $u\in\mathscr{E}_n$,  the solution $U(x,t;u,\overline U)$ of \eqref{A2} satisfies $U(\cdot,t;u,\overline U)\in\mathscr{E}_n$ for all $t\in[0,+\infty)$.
\end{proposition}
\proof Denote
\begin{align}\label{LU}
L(U):=\gamma(V)U''+\left(2\gamma'(V)V'+c\right)U'+\left(\gamma''(V)(V')^2+\gamma'(V)(V-U-cV')+a\right)U-bU^2
\end{align}
with $V$ defined in \eqref{2.23}. Noticing $\gamma(V)>0$, we have
$$U''+F(U,U')=\frac{L(U)}{\gamma(V)}.$$
Hence a function $U(x)$ is a super-solution (resp. sub-solution) of \eqref{A2} if $L(U)\leq0$ (reps. $L(U)\geq0$). Firstly we need to prove that for any solution $u\in\mathscr{E}_n$,
there exists $U(x,t;u,\overline U)\leq\overline{U}$.
For any $s\geq0$, from the definition of $\gamma(\cdot)$, we have
\begin{align}
0<\gamma(s)=\frac{1}{(1+s)^m}\leq1,\quad -m<\gamma'(s)=-\frac{m}{(1+s)^{m+1}}<0,\label{1.9}
\end{align}
and
\begin{align}
0<\gamma''(s)=\frac{m(m+1)}{(1+s)^{m+2}}\leq m(m+1).\label{1.13}
\end{align}
From \eqref{LU}, using \eqref{2.1} \eqref{2.2}, \eqref{1.9} and \eqref{1.13}, it is easy to verify that
\begin{align}
L(\eta )&=\left(\gamma''(V)(V')^2+\gamma'(V)(V-\eta-cV')+a\right)\eta -b\eta ^2\nonumber
\\
&\leq\left(\frac{m(m+1)}{1+a}\eta ^2+m\left(1+\frac{c}{\sqrt{1+a}}\right)\eta +a-b\eta \right)\eta .\label{1.6}
\end{align}
For $c\in[2\sqrt{a},c^*(a,b,m)]$, it is easy to verify that
\begin{align}
b-m\left(1+\frac{c}{\sqrt{1+a}}\right)\geq2\sqrt{\frac{m(m+1)a}{1+a}}.\label{1.16}
\end{align}
By the definition of $\eta$ in \eqref{2.8},
%\begin{align}
%\eta =\frac{b-m(1+\frac{c}{\sqrt{1+a}})-\sqrt{(b-m(1+\frac{c}{\sqrt{1+a}}))^2-\frac{4m(m+1)a}{1+a}}}{\frac{2m(m+1)}{1+a}},\label{2.8}
%\end{align}
one can check that
\begin{align}
\frac{m(m+1)}{1+a}\eta ^2+m\eta \left(1+\frac c{\sqrt{1+a}}\right)+a-b\eta =0.\label{1.15}
\end{align}
From \eqref{1.6} and \eqref{1.15}, we obtain $L(\eta )\leq0$.
On the other hand,
from \eqref{LU}, using \eqref{2.1}, \eqref{2.2},  \eqref{1.9}, and \eqref{1.13}, we obtain
\begin{eqnarray}\label{1.5}
\begin{aligned}
&L(e^{-\lambda x})\\
=&\gamma(V)\lambda^2e^{-\lambda x}-\left(2\gamma'(V)V'+c\right)\lambda e^{-\lambda x}+\left(\gamma''(V)(V')^2e^{-\lambda x}+\gamma'(V)(V-e^{-\lambda x}-cV')\right)e^{-\lambda x}\\
&\quad+ae^{-\lambda x}-be^{-2\lambda x}
\\
\leq&\lambda^2e^{-\lambda x}+\frac{2m\lambda}{\sqrt{1+a}}e^{-2\lambda x}-c\lambda e^{-\lambda x}+\frac{m(m+1)\eta }{1+a}e^{-2\lambda x}\quad+\left(m+\frac{cm}{\sqrt{1+a}}\right)e^{-2\lambda x}
\\
&\quad+ae^{-\lambda x}-be^{-2\lambda x}
\\
=&(\lambda^2-c\lambda+a)e^{-\lambda x}+\left(\frac{2m\lambda}{\sqrt{1+a}}+\frac{m(m+1)\eta }{1+a}+\frac{cm}{\sqrt{1+a}}+m-b\right)e^{-2\lambda x}.
\end{aligned}
\end{eqnarray}
From  \eqref{2.8} and \eqref{1.16}, it follows that
\begin{align}\label{2.11}
\eta &=\frac{2a}{b-m(1+\frac{c}{\sqrt{1+a}})+\sqrt{(b-m(1+\frac{c}{\sqrt{1+a}}))^2-\frac{4m(m+1)a}{1+a}}}
\leq\sqrt{\frac{a(1+a)}{m(m+1)}},
\end{align}
which along with the fact $c\in[2\sqrt{a},c^*(a,b,m)]$
implies
\begin{align*}
\frac{2m\lambda}{\sqrt{1+a}}+\frac{m(m+1)\eta }{1+a}+\frac{cm}{\sqrt{1+a}}+m-b
\leq\frac{(2\sqrt{a}+c)m}{\sqrt{1+a}}+\sqrt{\frac{am(m+1)}{1+a}}+m-b\leq0.
\end{align*}
Then from \eqref{1.5}, we obtain $L(e^{-\lambda x})\leq0$.
By the comparison principle for parabolic equations, it follows that $U(x,t;u,\overline U)\leq \overline{U}$.

Now we  prove that for any $u\in\mathscr{E}_n$, we have $U(x,t;u,\overline U)\geq \underline{U}_n$.
From \eqref{LU}, using \eqref{2.1}, \eqref{2.2}, \eqref{1.9} and \eqref{1.13}, we obtain
\begin{align}
L(\delta)&=\gamma''(V)(V')^2\delta+\gamma'(V)(V-\delta-cV')\delta
+\delta(a-b\delta)\label{2.27}
\\
&\geq\gamma'(V)(V-cV')\delta+\delta(a-b\delta)\nonumber
\\
&>\delta\left(a-b\delta-m\eta \left(1+\frac c{\sqrt{1+a}}\right)\right).\nonumber
\end{align}
Owing to the fact $c\in[2\sqrt{a},c^*(a,b,m)]$,
we have $b\geq3m\left(1+\frac c{\sqrt{1+a}}\right)$ and
\begin{align}
\eta <\frac{2a}{b-m\left(1+\frac c{\sqrt{1+a}}\right)}\leq\frac{a}{m\left(1+\frac c{\sqrt{1+a}}\right)}.\label{2.26}
\end{align}
Substituting \eqref{2.26} into \eqref{2.27}, we have $L(\delta)\geq0$ for sufficiently small $\delta$.
On the other hand, using \eqref{LU}, by direct calculations, we have
\begin{eqnarray}\label{1.34}
\begin{aligned}
&\quad L(d_ne^{-\theta_1(x)x}+d_0e^{-\theta_2(x)x})
\\
&=\gamma(V)(d_ne^{-\theta_1(x)x}+d_0e^{-\theta_2(x)x})''
+(2\gamma'(V)V'+c)(d_ne^{-\theta_1(x)x}+d_0e^{-\theta_2(x)x})'
\\
&\quad+\left(\gamma''(V)(V')^2+\gamma'(V)(V-(d_ne^{-\theta_1(x)x}+d_0e^{-\theta_2(x)x})-cV')+a\right)(d_ne^{-\theta_1(x)x}+d_0e^{-\theta_2(x) x})\\
&\quad-b(d_ne^{-\theta_1(x)x}+d_0e^{-\theta_2(x)x}))^2
\\
&\geq\gamma(V)d_ne^{-\theta_1(x)x}\left((\theta_1'(x)x)^2+\theta_1^2(x)+2\theta_1'(x)\theta_1(x)x-\theta_1''(x)x-2\theta_1'(x)\right)
\\
&\quad+\gamma(V)d_0e^{-\theta_2(x)x}\left((\theta_2'(x)x)^2+\theta_2^2(x)+2\theta_2'(x)\theta_2(x)x-\theta_2''(x)x-2\theta_2'(x)\right)
\\
&\quad+(2\gamma'(V)V'+c)\left(d_ne^{-\theta_1(x)x}(-\theta_1'(x)x-\theta_1(x))+d_0e^{-\theta_2(x)x}(-\theta_2'(x)x-\theta_2(x))\right)
\\
&\quad+\left(\gamma''(V)(V')^2+\gamma'(V)(V-cV')+a\right)(d_ne^{-\theta_1(x)x}+d_0e^{-\theta_2(x)x})\\
&\quad-b(d_ne^{-\theta_1(x)x}+d_0e^{-\theta_2(x)x})^2
\\
&=\left(\gamma(V)\theta_1^2(x)-c\theta_1(x)+a\right)d_ne^{-\theta_1(x)x}
+\left(\gamma(V)\theta_2^2(x)-c\theta_2(x)+a\right)d_0e^{-\theta_2(x)x}
\\
&\quad+d_ne^{-\theta_1(x)x}\left[\gamma(V)\left((\theta_1'(x)x)^2+2\theta_1'(x)\theta_1(x)x-\theta_1''(x)x-2\theta_1'(x)\right)-c\theta_1'(x)x\right.
\nonumber
\\
&\quad\left.-2\gamma'(V)V'(\theta_1'(x)x+\theta_1(x))+\gamma''(V)(V')^2+\gamma'(V)(V-cV')\right]\nonumber
\\
&\quad+d_0e^{-\theta_2(x)x}\left[\gamma(V)\left((\theta_2'(x)x)^2+2\theta_2'(x)\theta_2(x)x-\theta_2''(x)x-2\theta_2'(x)\right)-c\theta_2'(x)x\right.\nonumber
\\
&\quad\left.-2\gamma'(V)V'(\theta_2'(x)x+\theta_2(x))+\gamma''(V)(V')^2+\gamma'(V)(V-cV')\right]\\
&\quad-b(d_ne^{-\theta_1(x)x}+d_0e^{-\theta_2(x)x})^2.\nonumber
\end{aligned}
\end{eqnarray}
To prove that $L(d_ne^{-\theta_1(x)x}+d_0e^{-\theta_2(x)x})\geq0$, we consider the cases $c=2\sqrt{a}$
and $c>2\sqrt{a}$ separately.

\textbf{Case 1. $c=2\sqrt{a}$.}
In this case we have $d_0=1$ and substitute it into \eqref{1.34}. Using \eqref{1.9}, Lemma \ref{Lemma2.1} and Lemma \ref{Lemma2.3}, by choosing sufficiently small $\delta$, for $x>x_\delta$, we obtain
\begin{align*}
\gamma'(V)<0,\quad\gamma''(V)>0,\quad\theta_1'(x)>0,\quad\theta_1''(x)<0
\end{align*}
and
 \begin{align*}
\gamma(V)\theta_1^2(x)-c\theta_1(x)+a\geq0,\quad
\gamma(V)\theta_2^2(x)-c\theta_2(x)+a\geq \frac{a}{64},
\end{align*}
from which we obtain that for any $x>x_\delta$,
\begin{align}
&\quad L(d_ne^{-\theta_1(x)x}+e^{-\theta_2(x)x})\label{1.39}
\\
&\geq \frac{a}{64}e^{-\theta_2(x)x}+d_ne^{-\theta_1(x)x}\left[-2\gamma(V)\theta_1'(x)-c\theta_1'(x)x-2\gamma'(V)V'(\theta_1'(x)x+\theta_1(x))+\gamma'(V)(V-cV')\right]\nonumber
\\
&\quad+e^{-\theta_2(x)x}\left[-2\gamma(V)\theta_2'(x)-c\theta_2'(x)x-2\gamma'(V)V'(\theta_2'(x)x+\theta_2(x))+\gamma'(V)(V-cV')\right]\nonumber
\\
&\quad-b(d_ne^{-\theta_1(x)x}+e^{-\theta_2(x)x})^2.\nonumber
\end{align}
Furthermore,
from \eqref{1.22}, Lemma \ref{Lemma2.1} and Lemma \ref{Lemma2.2}, we have
\begin{align*}
0<\theta_1(x)<\sqrt{a},\quad 0<\theta_1'(x)\leq2K_1e^{-\frac\lambda2 x}
\end{align*}
and
\begin{align*}
0<V(x;u)\leq\min\left\{\frac{e^{-\lambda x}}{1+a},\eta \right\},\quad
|V'(x;u)|\leq \min\left\{\frac{\eta }{\sqrt{1+a}},\frac{e^{-\lambda x}}{\sqrt{1+a}}\right\}.
\end{align*}
By the above estimates and \eqref{1.9}, we arrive at the following estimates
\begin{align}
&-2\gamma(V)\theta_1'(x)-c\theta_1'(x)x-2\gamma'(V)V'(\theta_1'(x)x+\theta_1(x))+\gamma'(V)(V-cV')\label{1.35}
\\
\geq&-\left(4+2cx+\frac{4m\eta x}{\sqrt{1+a}}\right)K_1e^{-\frac{\lambda}{2}x}-\left(\frac{2m\sqrt{a}}{\sqrt{1+a}}+\frac{m}{1+a} +\frac{cm}{\sqrt{1+a}}\right)e^{-\lambda x}\nonumber
\\
=&-\left(4+4\sqrt{a}x+\frac{4m\eta x}{\sqrt{1+a}}\right)K_1e^{-\frac{\lambda}{2} x}-\left(\frac{4m\sqrt{a}}{\sqrt{1+a}}+\frac{m}{1+a}\right)e^{-\lambda x}. \nonumber
\end{align}
Then from the fact that $\theta_2(x)=\theta_1(x)+\frac\lambda 4$, we get
\begin{align}
&-2\gamma(V)\theta_2'(x)-c\theta_2'(x)x-2\gamma'(V)V'(\theta_2'(x)x+\theta_2(x))+\gamma'(V)(V-cV')\label{1.36}
\\
=&-2\gamma(V)\theta_1'(x)-c\theta_1'(x)x-2\gamma'(V)V'(\theta_1'(x)x+\theta_1(x))+\gamma'(V)(V-cV')-\frac12\gamma'(V)V'\lambda\nonumber
\\
\geq&-\left(4+4\sqrt{a}x+\frac{4m\eta x}{\sqrt{1+a}}\right)K_1e^{-\frac{\lambda}{2}  x}-\left(\frac{4m\sqrt{a}}{\sqrt{1+a}}+\frac{m}{1+a}+\frac{m\sqrt{a}}{2\sqrt{1+a}}\right)e^{-\lambda x}.\nonumber
\end{align}
Substituting \eqref{1.35} and \eqref{1.36} into \eqref{1.39}, we end up with
\begin{align}
&L(d_ne^{-\theta_1(x)x}+e^{-\theta_2(x)x})\nonumber
\\
&\geq e^{-\theta_2(x)x}\bigg\{\frac{a}{64}-K_1\left(4+4\sqrt{a}x+\frac{4m\eta x}{\sqrt{1+a}}\right)\left(e^{-\frac{\lambda}{2} x}+d_ne^{(\theta_2(x)-\theta_1(x)-\frac{\lambda}{2})x}\right)\label{1.51}
\\
&-\left(\frac{4m\sqrt{a}}{\sqrt{1+a}}+\frac{m}{1+a}\right)\left(e^{-\lambda x}+d_ne^{(\theta_2(x)-\theta_1(x)-\lambda)x}\right)-\frac{m\sqrt{a}}{2\sqrt{1+a}}e^{-\lambda x}\nonumber
\\
&-b\left(d_n^2e^{(\theta_2(x)-2\theta_1(x))x}+e^{-\theta_2(x)x}+2d_ne^{-\theta_1(x)x}\right)\bigg\}.\nonumber
\end{align}
From \eqref{1.24} and $\eqref{1.30}$, we have $\theta_2(x)-2\theta_1(x)<0$ and $\theta_2(x)-\theta_1(x)-\lambda<\theta_2(x)-\theta_1(x)-\frac{\lambda}{2}<0$ for $x>x_\delta$, then for $c=2\sqrt{a}$, by choosing $\delta$ sufficiently small in \eqref{1.51}, we obtain
$$L(d_ne^{-\theta_1(x)x}+e^{-\theta_2(x) x})\geq0$$ for all $x>x_\delta$.

\textbf{Case 2. $c>2\sqrt{a}$.}
Inserting $d_0=-1$ in \eqref{1.34}, using Lemma \ref{Lemma2.1}, Lemma \ref{Lemma2.2} and Lemma \ref{Lemma2.3}, we obtain

\begin{align*}
0<\theta_1(x)<\sqrt{a},\quad 0<\theta_1'(x)\leq2K_2e^{-\lambda x}\qquad 0>\theta_1''(x)\geq-2\lambda K_2e^{-\lambda x},
\end{align*}
\begin{align*}
0<V(x;u)\leq\min\left\{\frac{e^{-\lambda x}}{1+a},\eta \right\},\quad
|V'(x;u)|\leq \min\left\{\frac{\eta }{\sqrt{1+a}},\frac{e^{-\lambda x}}{\sqrt{1+a}}\right\}
\end{align*}
and
\begin{align*}
\gamma(V)\theta_1^2(x)-c\theta_1(x)+a\geq0,\quad
\gamma(V)\theta_2^2(x)-c\theta_2(x)+a\leq-\frac{\lambda(c-2\lambda)}{4k_0},
\end{align*}
from which and \eqref{1.9}, \eqref{1.13}, for any $x>x_\delta$, noticing that $\theta_2(x)=\theta_1(x)+\frac{\lambda}{k_0}$, we obtain
\begin{align}
&\quad L(d_ne^{-\theta_1(x)x}-e^{-\theta_2(x)x})\nonumber
\\
&\geq \left(\gamma(V)\theta_1^2(x)-c\theta_1(x)+a\right)d_ne^{-\theta_1(x)x}
-\left(\gamma(V)\theta_2^2(x)-c\theta_2(x)+a\right)e^{-\theta_2(x)x}\nonumber
\\
&\quad+d_ne^{-\theta_1(x)x}\left[-2\gamma(V)\theta_1'(x)-c\theta_1'(x)x-2\gamma'(V)V'(\theta_1'(x)x+\theta_1(x))+\gamma'(V)(V-cV')\right]\nonumber
\\
&\quad-e^{-\theta_2(x)x}\left[\gamma(V)((\theta_2'(x)x)^2+2\theta_2'(x)\theta_2(x)x-\theta_2''(x)x)-2\gamma'(V)V'(\theta_2'(x)x+\theta_2(x))\right.\nonumber
\\
&\quad\left.+\gamma''(V)(V')^2-c\gamma'(V)V'\right]-b(d_ne^{-\theta_1(x)x}-e^{-\theta_2(x)x})^2\label{1.53}
\\
&\geq e^{-\theta_2(x)x}\left\{\frac{\lambda(c-2\lambda)}{4k_0}-b\left(d_n^2e^{(\theta_2(x)-2\theta_1(x))x}+e^{-\theta_2(x)x}\right)\right.\nonumber
\\
&\quad-\left(4K_2+2cxK_2+\frac{2m}{\sqrt{1+a}}(2K_2e^{-\lambda x}x+\sqrt{a})+m\left(\frac{1}{1+a}+\frac{c}{\sqrt{1+a}}\right)\right) d_ne^{(\theta_2(x)-\theta_1(x)-\lambda)x}\nonumber
\\
&\quad \left.-\left((2K_2x)^2e^{-\lambda x}+4K_2\left(\sqrt{a}+\frac{\lambda}{k_0}\right)x+2\lambda K_2x+\frac{2m}{\sqrt{1+a}}\left(2K_2xe^{-\lambda x}+\sqrt{a}+\frac{\lambda}{k_0}\right)\right.\right.\nonumber
\\
&\qquad\left.\left.+\frac{m(m+1)}{1+a}e^{-\lambda x}+\frac{cm}{\sqrt{1+a}}\right)e^{-\lambda x}\right\}.\nonumber
\end{align}
Noticing $\theta_2(x)-2\theta_1(x)<0$ and $\theta_2(x)-\theta_1(x)-\lambda<0$ for $x>x_\delta$, then for $c>2\sqrt{a}$, by choosing $\delta$ sufficiently small in \eqref{1.53}, we obtain
$$L(d_ne^{-\theta_1(x)x}-e^{-\theta_2(x) x})\geq0$$ for all $x>x_\delta$.
Then by the comparison principle for parabolic equations, we obtain $U(x,t;u)\geq \underline{U}_n$ for $c\geq2\sqrt{a}$.

Summing up, by choosing
$$\eta:=\frac{b-m(1+\frac{c}{\sqrt{1+a}})-\sqrt{(b-m(1+\frac{c}{\sqrt{1+a}}))^2-\frac{4m(m+1)a}{1+a}}}{\frac{2m(m+1)}{1+a}}$$
and sufficiently small $\delta$, $(\overline{U},\underline{U}_n)$ is a pair of super- and sub-solutions of \eqref{A2} (see a schematic of super- and sub-solutions illustrated in Fig.\ref{fig}).
Denoting $U(x,t;u,\overline U)$ the unique solution of \eqref{A2}, by the comparison principle for parabolic equations, we obtain $\underline{U}_n\leq U(x,t;u,\overline U)\leq \overline{U}$ and thus $U(x,t;u,\overline U)\in\mathscr{E}_n$. This completes the proof of Lemma \ref{Pro2.2}.

$\hfill\square$

%\begin{figure}[htbp]
%\centering
%\includegraphics[width=7cm,angle=0]{super-sub}
%\caption{A schematic of super- and sub-solutions, where the solid black line represents the super-solution $\overline{U}(x)$ and the dashed red line is the sub-solution $\underline{U}_n(x)$.}
%\label{fig}
%\end{figure}

\subsection{An auxiliary elliptic problem}
Now for $u\in\mathscr{X}_0:=\bigcap_{n>1}\mathscr{E}_n$, we study the following problem
\begin{align}\label{1.4}
\left\{
\begin{array}{ll}
\gamma(V)U''+\left(2\gamma'(V)V'+c\right)U'+\left(\gamma''(V)(V')^2+\gamma'(V)(V-U-cV')+a\right)U-bU^2=0, \\
V''+cV'+u-V=0,
\end{array}
\right.
\end{align}
which  is equivalent to solving $L(U)=0$.

\begin{proposition}\label{proposition2.1}
For every $u\in\mathscr{X}_0$, $b\geq b^*(m,a)$ and $c\in[2\sqrt{a},c^*(a,b,m)]$, denote $U(x,t;u,\overline{U})$ the solution of \eqref{A2} with $U(x,0;u,\overline{U})=\overline{U}$, there exists a unique function $U(x;u)\in\mathscr{X}_0$ such that
$$U(x;u)=\lim_{t\to\infty}U(x,t;u,\overline{U})=\inf_{t>0}U(x,t;u,\overline{U})$$
and $U(x;u)$ is the unique solution of \eqref{1.4} satisfying
\begin{align}
\liminf_{x\to-\infty}U(x;u)>0\quad\hbox{and}\;\lim_{x\to+\infty}\frac{U(x;u)}{e^{-\lambda x}}=1.\label{1.37}
\end{align}
\end{proposition}

\proof From Proposition \ref{Pro2.2}, we have
\begin{align}
U(x,t;u,\overline{U})\leq\overline{U}(x)\quad\hbox{for all}\;(x,t)\in\mathbb R\times[0,+\infty).\label{1.10}
\end{align}
For any $0\leq t_1\leq t_2$, noticing
$$U(x,t_2;u,\overline{U})=U(x,t_1;u,U(x,t_2-t_1;u,\overline{U})),$$ from \eqref{1.10}, we have
$$U(x,t_2-t_1;u,\overline{U})\leq\overline{U}(x),$$ then using again the comparison principle for parabolic equations, we obtain
$$U(x,t_2;u,\overline{U})\leq U(x,t_1;u,\overline{U}),$$
which implies that $U(x,\cdot;u,\overline{U})$ is decreasing with respect to $t$. Noticing that $U(x,\cdot;u,\overline{U})$ has lower and upper bounds since $U(x,\cdot;u,\overline{U})\in \mathscr{E}_n$ as shown in Lemma \ref{Pro2.2}, one can conclude that there exists a unique $U(x;u)$ such that
\begin{align}
U(x;u)=\lim_{t\to\infty}U(x,t;u,\overline{U})=\inf_{t>0}U(x,t;u,\overline{U})
\label{1.11}
\end{align}
for all $x\in\mathbb R$. Denote $$U_n(x,t)=U(x,t+t_n;u,\overline{U})$$ for $(x,t)\in\mathbb R\times[0,\infty)$, where $\{t_n\}_{n\geq1}$ is an increasing
sequence of positive real numbers converging to $+\infty$.
Then from the elliptic regularity theory for \eqref{A1} and parabolic regularity theory for \eqref{A2} (cf. \cite{Lady}), we obtain that for all $1<p<\infty$, $R>0$, $T>0$,
$$\|V\|_{W^{2,p}(-R,R)}\leq C$$
and
$$\|U_n\|_{W^{2,1}_{p}((-R,R)\times(0,T))}\leq C.$$
From Sobolev embedding theorem, we obtain
$$\|V\|_{C_{\mathrm{loc}}^{1,\alpha}(\R)}\leq C$$
and
$$\|U_n\|_{C^{\alpha,\alpha/2}_{\mathrm{loc}}(\mathbb R\times(0,+\infty))}\leq C.$$
The Arzel\`{a}-Ascoli's theorem and Schauder's theory for parabolic equation imply that there is a subsequence $\{U_{n'}\}_{n'\geq1}$ of the sequence $\{U_{n}\}_{n\geq1}$ and a function $\tilde{U}\in C^{2,1}(\mathbb R\times(0,\infty))$, such that $\{U_{n'}\}_{n'\geq1}$ converges to $\tilde{U}$ locally uniformly in $C^{2,1}(\mathbb R\times(0,\infty))$ as $n'\to\infty$. Hence $\tilde{U}(x,t)$ solves \eqref{1.4}
and $\tilde{U}\in\mathscr{X}_0$. On the other hand, noticing $\tilde{U}(x,t)=\lim_{t\to\infty}U(x,t;u,\overline{U})$, from \eqref{1.11}, we have $U(x;u)=\tilde{U}(x,t)$ for every $x\in\mathbb R$ and $t\geq0$, from which we obtain that $U(x;u)\in\mathscr{X}_0$ is a solution of \eqref{1.4}. Furthermore, from \eqref{1.47}
%that is
%$$\lim_{x\to+\infty}e^{(\theta_1(x)-\lambda)x}=1,$$
and the definition of $\mathscr{X}_0$, we obtain
\begin{align}
\liminf_{x\to-\infty}U(x;u)>0\label{1.44}
\end{align}
and
\begin{align}
d_n\leq\liminf_{x\to+\infty}\frac{U(x;u)}{e^{-\lambda x}}\leq\limsup_{x\to+\infty}\frac{U(x;u)}{e^{-\lambda x}}=1\label{1.42}
\end{align}
for any $n\geq2$. Noticing $\lim_{n\to\infty}d_n=1$, by taking $n\to\infty$ in \eqref{1.42}, we obtain
\begin{align}
\lim_{x\to+\infty}\frac{U(x;u)}{e^{-\lambda x}}=1.\label{1.43}
\end{align}
The uniqueness of $U(x;u)$ satisfying \eqref{1.37} follows from the same arguments as that in Lemma 3.6 in \cite{SS1}.
The proof is thus completed. $\hfill\square$

\section{Proof of main theorems}
In this section, we shall prove Theorem \ref{main-thm1} and Theorem \ref{main-thm2}. To this end, we first prove the following result concerning the asymptotic behavior of solutions to \eqref{1.2} as $z \to \pm \infty$.
%In preceding subsection, we establish the existence of solutions of \eqref{1.4} for given $u\in \mathscr{X}_0$. Below we shall proceed to study the asymptotic behavior of solutions to \eqref{1.4} as $z\to \pm\infty$.
\begin{proposition}\label{Pro2.3}
Assume that $a>0$ and $m>0$ satisfy \eqref{1.12}.
%\begin{equation}\label{1.12}
%m\sqrt{\frac{a(1+a)}{m(m+1)}}\left(\sqrt{\frac{a(1+a)}{m(m+1)}}+1\right)^m<1.
%\end{equation}
Then any solution $(U,V)\in(C^2(\mathbb R)\cap \mathscr{X}_0)^2$ to \eqref{1.2} has the property that
$$\lim_{z\to+\infty}U(z)=\lim_{z\to+\infty}V(z)=0,\quad\lim_{z\to-\infty}U(z)
=\lim_{z\to-\infty}V(z)=a/b$$
and
$$\lim_{z\to\pm\infty}U'(z)=\lim_{z\to\pm\infty}V'(z)=0.$$
\end{proposition}

\proof From the fact that $(U,V)\in \mathscr{X}_0^2$ and Lemma \ref{Lemma2.2}, we obtain
\begin{align}
|U(z)|\leq \eta , \quad|V(z)|\leq \eta \quad \hbox{and}\; |V'(z)|\leq \frac{\eta }{\sqrt{1+a}}\label{1.54}
\end{align}
for all $z\in\mathbb R$. From the first equation of \eqref{1.2}, by the H\"older regularity estimates for bounded solutions of elliptic equations and the standard Schauder theory, there exists $C>0$ independent of $z$ and $\alpha\in(0,1)$ such that $\|U\|_{C^{2,\alpha}(z,z+1)}\leq C$ and $\|V\|_{C^{2,\alpha}(z,z+1)}\leq C$ for all $z\in\mathbb R$, from which it follows that
\begin{align}
|U'(z)|\leq C, \quad|U''(z)|\leq C\quad\hbox{and}\;|V''(z)|\leq C\label{1.49}
\end{align}
for all $z\in\mathbb R$. Multiplying the first equation of \eqref{1.2} by $(a-bU)$, integrating over $[-R,R]$, we obtain
\begin{align*}
0=&\int_{-R}^R(\gamma(V)U)''(a-bU)dz+c\int_{-R}^RU'(a-bU)dz+\int_{-R}^R U(a-bU)^2dz
\\
=&(\gamma(V)U)'(a-bU)\big|_{z=-R}^{z=R}+b\int_{-R}^R(\gamma'(V)V'U+\gamma(V)U')U'dz
+caU\big|_{z=-R}^{z=R}
\\
&-\frac12cbU^2\big|_{z=-R}^{z=R}+\int_{-R}^R U(a-bU)^2dz.
\end{align*}
Then using \eqref{1.9}, \eqref{1.54} and \eqref{1.49}, we find a constant $C_1$ independent of $R$ such that
\begin{align}
&\frac{b}{(1+\eta )^m}\int_{-R}^R|U'|^2dz+\int_{-R}^R U(a-bU)^2dz\nonumber
\\
\leq &b\int_{-R}^R\gamma(V)|U'|^2dz+\int_{-R}^R U(a-bU)^2dz\label{r}
\\
\leq &b\int_{-R}^R|\gamma'(V)V'UU'|dz-(\gamma(V)U)'(a-bU)\big|_{z=-R}^{z=R}-caU\big|_{z=-R}^{z=R}+\frac12cbU^2\big|_{z=-R}^{z=R}\nonumber
\\
\leq &C_1+\frac12bm\eta \left(\int_{-R}^R|U'|^2dz+\int_{-R}^R|V'|^2dz\right).\nonumber
\end{align}
On the other hand, multiplying the second equation of \eqref{1.2} by $V''$ and integrating the result over $[-R,R]$, we obtain
\begin{align*}
0&=\int_{-R}^R|V''|^2dz+c\int_{-R}^RV'V''dz+\int_{-R}^RUV''dz-\int_{-R}^RVV''dz
\\
&=\int_{-R}^R|V''|^2dz+\frac{c}{2}(V')^2\big|_{z=-R}^{z=R}+UV'\big|_{z=-R}^{z=R}
-\int_{-R}^RU'V'dz-VV'\big|_{z=-R}^{z=R}+\int_{-R}^R|V'|^2dz.
\end{align*}
This along with \eqref{1.54} and \eqref{1.49} yields
\begin{align*}
\int_{-R}^R|V''|^2dz+\int_{-R}^R|V'|^2dz&\leq C_2+\int_{-R}^RU'V'dz
\leq C_2+\frac12\int_{-R}^R|U'|^2dz+\frac12\int_{-R}^R|V'|^2dz,
\end{align*}
where $C_2$ is  a constant independent of $R$. Then it follows that
\begin{align}
\int_{-R}^R|V'|^2dz\leq 2C_2+\int_{-R}^R|U'|^2dz.\label{rr}
\end{align}
Substituting \eqref{rr} into \eqref{r}, one can find  a constant $C_3=C_1+bm\eta C_2$ independent of $R$ such that
\begin{align}
\frac{b}{(1+\eta )^m}\int_{-R}^R|U'|^2dz+\int_{-R}^R U(a-bU)^2dz
\leq C_3+bm\eta \int_{-R}^R|U'|^2dz.\label{2.9}
\end{align}
Note that \eqref{2.11} together with condition \eqref{1.12} implies
\begin{equation}
\frac{1}{(1+\eta )^m}-m\eta>0.
\end{equation}
Sending $R\to\infty$ in \eqref{2.9}, we obtain
\begin{align}
b\left(\frac{1}{(1+\eta )^m}-m\eta \right)\int_{\mathbb R}|U'|^2dz+\int_{\mathbb R} U(a-bU)^2dz
\leq C_3.\label{cc1}
\end{align}
By sending $R\to\infty$ in \eqref{rr}, we find a constant $C_4>0$ such that
\begin{align}\label{rrr}
\int_{\mathbb R}|V'|^2dz\leq C_4.
\end{align}
Then \eqref{cc1} and \eqref{rrr} assert that
\begin{align}
U'\in L^2(\mathbb R),\ \ U(a-bU)^2\in L^1(\mathbb R), \ \ V'\in L^2(\mathbb R).\label{2.12}
\end{align}
From \eqref{1.49} and \eqref{2.12} , we obtain
\begin{align}\label{sss}
\lim_{z\to\pm\infty}U(z)\in\{0,a/b\},\quad \lim_{z\to\pm\infty}U'(z)=0\quad\hbox{and}\;
\lim_{z\to\pm\infty}V'(z)=0.
 \end{align}
Furthermore, from the definition of $\mathscr{X}_0$ and the fact that $U\in\mathscr{X}_0$, we obtain $$\lim_{z\to+\infty}U(z)=0\quad\hbox{and}\; \lim_{z\to-\infty}U(z)=a/b.$$
On the other hand, from the second equation of \eqref{1.2}, we have
\begin{align}
V(z)=\frac{1}{\lambda_2-\lambda_1}\left(\int_{-\infty}^ze^{\lambda_1(z-s)}U(s)ds+\int_z^{+\infty}
e^{\lambda_2(z-s)}U(s)ds\right)\label{1.40}
\end{align}
with $\lambda_1<0$ and $\lambda_2>0$ defined in \eqref{1.38}.
Applying L'Hopital's rule to \eqref{1.40}, from the fact \eqref{sss}, we obtain
\begin{align*}
\lim_{z\to+\infty}V(z)&=\lim_{z\to+\infty}\frac{1}{\lambda_2-\lambda_1}
\left(\frac{\int_{-\infty}^ze^{-\lambda_1s}U(s)ds}{e^{-\lambda_1z}}
+\frac{\int_z^{+\infty}e^{-\lambda_2s}U(s)ds}{e^{-\lambda_2z}}\right)
\\
&=\frac{1}{\lambda_2-\lambda_1}\lim_{z\to+\infty}\left(\frac{U(z)}{-\lambda_1}+\frac{U(z)}{\lambda_2}\right)
\\
&=\lim_{z\to+\infty}U(z)=0
\end{align*}
and
\begin{align*}
\lim_{z\to-\infty}V(z)&=\lim_{z\to-\infty}\frac{1}{\lambda_2-\lambda_1}\left(\frac{\int_{-\infty}^ze^{-\lambda_1s}U(s)ds}{e^{-\lambda_1z}}
+\frac{\int_z^{+\infty}e^{-\lambda_2s}U(s)ds}{e^{-\lambda_2z}}\right)
\\
&=\frac{1}{\lambda_2-\lambda_1}\lim_{z\to-\infty}\left(\frac{U(z)}{-\lambda_1}+\frac{U(z)}{\lambda_2}\right)
\\
&=\lim_{z\to-\infty}U(z)=\frac ab.
\end{align*}
This completes the proof.

$\hfill\square$

\subsection{Proof of Theorem \ref{main-thm1}}
Note that a fixed point of the mapping $u\ni\mathscr{X}_0\mapsto U(\cdot,u)\in\mathscr{X}_0$ formed in \eqref{1.4} is a solution to the wave equations \eqref{1.2}. Hence to prove the existence of travelling wave solutions to \eqref{1.1}, it suffices to prove that the mapping $u\ni\mathscr{X}_0\mapsto U(\cdot,u)\in\mathscr{X}_0$ formed in \eqref{1.4} has a fixed point. We shall achieve this by the Schauder fixed point theorem.

First, we prove that the mapping $u\ni\mathscr{X}_0\mapsto U(\cdot,u)\in\mathscr{X}_0$ is compact. Let $\{u_n\}_{n\geq1}$ be a sequence in $\mathscr{X}_0$. Denote $U_n=U(\cdot,u_n)$, we have $U_n\in\mathscr{X}_0$.
From the elliptic regularity theorem,  we have that for all $p>1$
$$\|U_n\|_{W^{2,p}_{\mathrm{loc}}(\mathbb R)}\leq C.$$
From Sobolev embedding theorem, we obtain
$$\|U_n\|_{C^{\alpha}_{\mathrm{loc}}(\mathbb R)}\leq C,$$
the Arzela-Ascoli's theorem implies that there is a subsequence $\{U_{n'}\}_{n'\geq1}$ of the sequence $\{U_{n}\}_{n\geq1}$ and a function $U(x)\in C(\mathbb R)$, such that $\{U_{n'}\}_{n'\geq1} \to U(x)$  locally uniformly in $ C(\mathbb R)$. Furthermore, we have $U(x)\in\mathscr{X}_0$. Then  the mapping $u\ni\mathscr{X}_0\mapsto U(\cdot,u)\in\mathscr{X}_0$ is compact.

Second, we prove that the mapping $u\ni\mathscr{X}_0\mapsto U(\cdot;u)\in\mathscr{X}_0$ is continuous. To this end, denote
$$\|u\|_*=\sum_{n=1}^\infty\frac1{2^n}\|u\|_{L^\infty([-n,n])}.$$
Then any sequence of functions in $\mathscr{X}_0$ is convergent with respect to norm $\|\cdot\|_*$ if and only if it converges locally uniformly on $\mathbb R$. Let $u\in\mathscr{X}_0$ and
$\{u_n\}_{n\geq1}$ be a sequence in $\mathscr{X}_0$ such that $u_n$ converges to $u$ locally uniformly on $\mathbb R$ as $n\to\infty$. Then by the elliptic regularity theorem applied to the second equation of \eqref{1.4} and Sobolev embedding theorem, we obtain
$$\|V(\cdot;u_n)\|_{C^{1,\alpha}_{\mathrm{loc}}(\mathbb R)}\leq C.$$
Form the Arzel\`{a}-Ascoli's theorem, there exists a subsequence of $\{V(\cdot;u_{n})\}_{n\geq1}$, still denoted by itself without confusion, such that
$$\lim_{n'\to\infty}V(\cdot;u_{n})=V(\cdot;u)\quad\hbox{in}\; C^1_{\mathrm{loc}}(\mathbb R).$$
Suppose by contradiction that that the mapping $u\ni\mathscr{X}_0\mapsto U(\cdot;u)\in\mathscr{X}_0$ is not continuous, then there exists $\delta>0$ and a
subsequence $\{u_{n'}\}_{n'\geq1}$ such that
\begin{align}
\|U(\cdot; u_{n'})-U(\cdot; u)\|_*\geq\delta,\quad\forall n\geq1.\label{1.41}
\end{align}
By Schauder's theory applied to the first equation of \eqref{1.4} and Sobolev embedding theorem, from the Arzel\`{a}-Ascoli's theorem, there is a subsequence $\{U(\cdot; u_{n''})\}_{n''\geq1}$ of the sequence $\{U(\cdot;u_{n'})\}_{n'\geq1}$ and a function $U(\cdot)\in C^{2}(\mathbb R)$, such that $\{U(\cdot; u_{n''})\}_{n''\geq1}$ converges to $U(\cdot)$ in $C^2_{\mathrm{loc}}(\mathbb R)$ and $U$ is a solution of \eqref{1.4}. Moreover, from the fact that $U(\cdot;u_{n''})\in\mathscr{X}_0$ and
$$\lim_{n\to\infty}\|U(\cdot;u_{n''})- U(\cdot)\|_*=0,$$
we obtain
$U(\cdot)\in\mathscr{X}_0$.
Then from Proposition
\ref{proposition2.1}, we obtain $U(\cdot)=U(\cdot,u).$
By \eqref{1.41}, then
$$\|U(\cdot;u)-U(\cdot)\|_*\geq\delta,$$
which is a contradiction. Hence the mapping $u\ni\mathscr{X}_0\mapsto U(\cdot;u)\in\mathscr{X}_0$ is continuous.

Now by the Schauder's fixed point theorem, there is $U\in\mathscr{X}_0$ such that $U(\cdot)=U(\cdot;U)$. Denote $V(\cdot):=V(\cdot;U)$. Then $(U,V)$ is a solution of \eqref{1.2}. From the definition of $\mathscr{X}_0$ and \eqref{1.47}, we obtain $$\lim_{z\to+\infty}\frac{U(z)}{e^{-\lambda z}}=1.$$
This along with \eqref{1.45}-\eqref{1.46} and L'H\^{o}pital's Rule yields
\begin{align*}
\lim_{z\to+\infty}\frac{V(z)}{e^{-\lambda z}}
&=\lim_{z\to+\infty}\frac{1}{\lambda_2-\lambda_1}\left(\frac{\int_{-\infty}^ze^{-\lambda_1s}U(s)ds}{e^{-(\lambda_1+\lambda)z}}
+\frac{\int_z^{+\infty}e^{-\lambda_2s}U(s)ds}{e^{-(\lambda_2+\lambda)z}}\right)
\\
&=\frac{1}{\lambda_2-\lambda_1}\lim_{z\to+\infty}\left(\frac{U(z)}{-(\lambda_1+\lambda)e^{-\lambda z}}-\frac{U(z)}{-(\lambda_2+\lambda)e^{-\lambda z}}\right)=\frac1{1+a}.
\end{align*}
Since $U\in\mathscr{X}_0$, it follows that $\liminf\limits_{z\to-\infty}U(z)>0$. On the other hand, noticing for $z<x_\delta$, $U(z)>\delta$ and then
\begin{align*}
V(z)&=\frac{1}{\lambda_2-\lambda_1}\left(\int_{-\infty}^ze^{\lambda_1(z-s)}U(s)ds
+\int_z^{+\infty}e^{\lambda_2(z-s)}U(s)ds\right)
\\
&\geq\frac{\delta}{\lambda_2-\lambda_1}\int_{-\infty}^ze^{\lambda_1(z-s)}ds
\\
&=\frac{\delta}{(\lambda_2-\lambda_1)(-\lambda_1)}>0,
\end{align*}
from which $\liminf\limits_{z\to-\infty}V(z)>0$ follows. Finally by the assumption \eqref{1.12} and Proposition \ref{Pro2.3}, we finish the proof of Theorem \ref{main-thm1}.
%$$m\sqrt{\frac{a(1+a)}{m(m+1)}}\left(\sqrt{\frac{a(1+a)}{m(m+1)}}+1\right)^m<1,$$
%we obtain
%$$\lim_{z\to-\infty}U(z)=\lim_{z\to-\infty}V(z)=a/b.$$ The proof of Theorem \ref{main-thm1} is completed.

$\hfill\square$

\subsection{Proof of Theorem \ref{main-thm2}}
Arguing by contradiction, for $c<2\sqrt{a}$, we suppose that there is a travelling wave solution $(u(x,t),v(x,t))=(U(x\cdot\xi-ct),V(x\cdot\xi-ct))$ of \eqref{1.1} connecting the constant solutions $(a/b,a/b)$ and $(0,0)$. Take a sequence $\{z_n\}$ with $z_n\to+\infty$, then $$\lim_{n\to+\infty}U(z_n)=\lim_{n\to+\infty}V(z_n)=\lim_{n\to+\infty}V'(z_n)=0.$$ Now we set
$$h_n(z)=\frac{U(z+z_n)}{U(z_n)},\quad U_n(z)=U(z+z_n),\quad V_n(z)=V(z+z_n).$$
As $U$ is bounded and satisfies \eqref{1.2}, the Harnack inequality implies that the shifted function $U_n(z)$, $V_n(z)$ and $V_n'(z)$ converge to zero locally uniformly in $z$ and the sequence $h_n$ is locally uniformly bounded and satisfies
\begin{align*}
\left\{
\begin{array}{ll}
\gamma''(V_n)(V'_n)^2h_n+\gamma'(V_n)(V_n-U_n-cV_n')h_n+2\gamma'(V_n)V_n'h_n'+\gamma(V_n)h_n''+ch_n'+h_n(a-bU_n)=0,\\
V_n''+U_n-V_n+cV_n'=0
\end{array}
\right.
\end{align*}
in $\mathbb R$. Thus up to a subsequence, the sequence $\{h_n\}_{n\geq 1}$ converges to a function $h$ that satisfies
\begin{align}
h''+ch'+ah=0\quad\hbox{in}\;\mathbb R.\label{2.5}
\end{align}
Moreover, $h$ is nonnegative and $h(0)=1$. Equation \eqref{2.5} admits such a solution if and only if $c\geq 2\sqrt{a}$, which leads to a contradiction. This denies our assumption and hence \eqref{1.1} admits no travelling wave solution connecting $(a/b,a/b)$ and $(0,0)$ with speed $c<2\sqrt{a}$.$\hfill\square$

\section{Selection of wave profiles}
By introducing some auxiliary problems and spatially inhomogeneous relaxed decay rates for super- and sub-solutions constructed, we manage to establish the existence of traveling wavefront solutions to the density-suppressed motility system \eqref{1.1} with decay motility function \eqref{motility}, where we find that there is a minimal wave speed coincident with the one for the cornerstone Fisher-KPP equation, and a maximum wave speed $c$ resulting from the nonlinear diffusion.
%As far as we know, this is the first result on the traveling wave solutions to the density-suppressed motility system \eqref{1.1} with continuous motility function. Though we consider a specialized motility function \eqref{motility}, our ideas and analysis can be extended to other kind of motility functions with direct modifications.
However, we are unable to characterize further properties of wave profiles such as monotonicity, stability and so on. In this section, we shall discuss the selection of possible wave profiles motivated by some argument in \cite{ou-yuan}.

\subsection{Trailing edge wave profiles}
In the spatially homogeneous situation, the system \eqref{1.1} has equilibria $(0,0)$ and $(a/b, a/b)$, which are unstable saddle and stable node respectively. This suggests that we should look for traveling wavefront solutions to \eqref{1.1} connecting $(a/b, a/b)$ to $(0,0)$ as we have done in the paper. Now we linearize the ODE system \eqref{1.2} at the origin $(0,0)$ and let  $U'=X, V'=Y$. Then we get the following linear system of $(U,X,V,Y)$
\begin{equation}\label{linear1}
\begin{pmatrix}
U'\\X'\\V'\\Y'
\end{pmatrix}
=
\begin{pmatrix}
0 & 1 & 0 & 0\\
-\frac{a}{\gamma(0)} & -\frac{c}{\gamma(0)} & 0 & 0\\
0 & 0 & 0 & 1\\
-1 & 0 & 1 & -c
\end{pmatrix}
\begin{pmatrix}
U\\X\\V\\Y
\end{pmatrix}.
\end{equation}
The eigenvalue $\lambda$ of the above coefficient matrix is
$$\Big(\lambda^2+\frac{c}{\gamma(0)}\lambda +\frac{a}{\gamma(0)}\Big) \Big(\lambda^2+c\lambda -1 \Big)=0.$$
To ensure there is a positive trajectory connecting the equilibria $(0,0)$ and $(a/b,a/b)$, we need to rule out the case that $(0,0)$ is a spiral, which amounts to require
\begin{equation}\label{minimal}
c\geq 2 \sqrt{\gamma(0) a}.
\end{equation}
With $\gamma(v)$ given in \eqref{motility}, $\gamma(0)=1$ and \eqref{minimal} is equivalent to $c\geq 2 \sqrt{a}$. This is well consistent with our results obtained in Theorem \ref{main-thm1} and Theorem \ref{main-thm2}.  Under the restriction \eqref{minimal}, it can be easily check that the origin $(0,0)$ is either a stable node or saddle point, which indicates that the traveling wave profile around the origin $(0,0)$ will not be oscillatory or periodic.

Next we linearize the system \eqref{1.2} at $(a/b,a/b)$ and arrive at the following linearized system
\begin{equation}\label{linear2}
\begin{pmatrix}
U'\\X'\\V'\\Y'
\end{pmatrix}
=
\begin{pmatrix}
0 & 1 & 0 & 0\\
\frac{a(b+\sigma_2)}{\sigma_1 b} & -\frac{c}{\sigma_1} & -\frac{a \sigma_2}{b \sigma_1} & \frac{a \sigma_2 c}{b \sigma_1}\\
0 & 0 & 0 & 1\\
-1 & 0 & 1 & -c
\end{pmatrix}
\begin{pmatrix}
U\\X\\V\\Y
\end{pmatrix}
\end{equation}
where $\sigma_1=\gamma(a/b), \sigma_2=\gamma'(a/b)$.
% and $\sigma=\frac{a(b+\sigma_2)}{\sigma_1 b}\frac{a}{\sigma_1}\big(1+\frac{\sigma_2}{b}\big)$.
By some tedious computation, we find that the eigenvalue $\lambda$ of the above coefficient matrix is determined by the following characteristic equation
\begin{equation}\label{eigen2}
\lambda^4+\Big(c+\frac{c}{\sigma_1}\Big)\lambda^3+\bigg(\frac{c^2}{\sigma_1}-\frac{a(b+\sigma_2)}{\sigma_1 b}-1\bigg)\lambda^2-\frac{(a+1)c}{\sigma_1}\lambda+\frac{a}{\sigma_1}=0.
\end{equation}
We suppose that there are periodic solutions near the positive equilibrium $(a/b,a/b)$, namely  the above characteristic equation has purely imaginary roots $\lambda= \pm \omega i$, where $\omega$ is a real number. Then the substitution of this ansatz into the equation \eqref{eigen2} immediately yields a necessary condition $c=0$, and consequently we get
\begin{equation}\label{omega}
\omega^4-\bigg(\frac{a(b+\sigma_2)}{\sigma_1 b}+1\bigg)\omega^2+\frac{a}{\sigma_1}=0.
\end{equation}
Notice that $\sigma_2=\gamma'(a/b)<0$. Then a necessary and sufficient condition warranting that the equation \eqref{omega} has a real root $\omega$ is
%\begin{equation}\label{conp}
%\sigma_2<-b\Big(\frac{\sigma_1}{a}+1\Big).
%\end{equation}
\begin{equation}\label{conp}
|\sigma_2|<\frac{b}{a}\sigma_1 \bigg(\sqrt{\frac{a}{\sigma_1}}-1\bigg)^2.
\end{equation}
That is, the linearized system \eqref{linear2} at the equilibrium $(a/b,a/b)$ will have periodic solutions if the condition \eqref{conp} is fulfilled. Thereof we anticipate that the non-monotone traveling wave solutions oscillating about the critical point $(a/b,a/b)$ may exist, but whether the condition \eqref{conp} is sufficient to guarantee that the nonlinear system \eqref{1.1} has similar oscillatory behavior around the equilibrium $(a/b,a/b)$ is very hard to determine and even to predict due to the complexity induced by the nonlinear diffusion and cross-diffusion in the system. Below we shall use numerical simulations to illustrate that indeed the condition \eqref{conp} plays a critical role for the nonlinear system in determining the monotonicity of wave profiles.

\begin{figure}[htbp]
\centering
\includegraphics[width=14cm,angle=0]{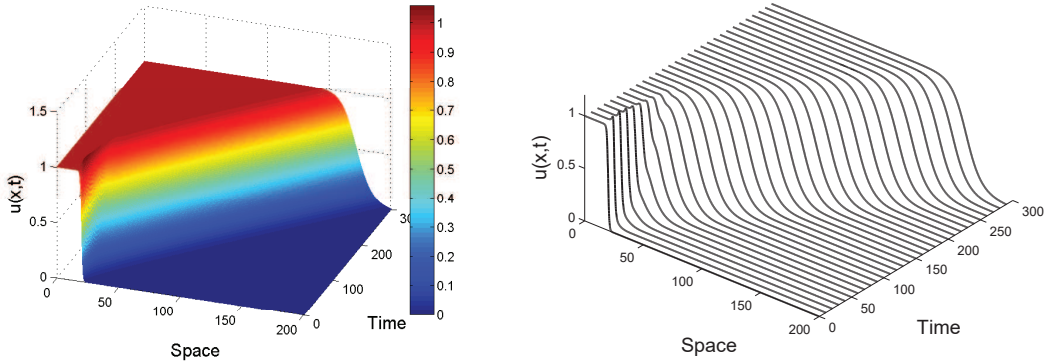}
\caption{Numerical simulations of wave propagation generated by the system \eqref{1.1} in $[0,200]$ with $\gamma(v)=\frac{1}{(1+v)^m}$ with $m=6$, $a=b=0.1, u_0=v_0=\frac{1}{1+e^{2(x-20)}}$.}
\label{fig1}
\end{figure}

We consider the motility function $\gamma(v)=\frac{1}{(1+v)^m} (m>0)$ as given in \eqref{1.4}. With simple calculation, we find that the condition \eqref{conp} amounts to
\begin{equation}\label{conn}
\sqrt{m}<\sqrt{\frac{1+\vartheta}{\vartheta}} \ \big| \sqrt{a(1+\vartheta)^m}-1\big|,\ \ \vartheta=\frac{a}{b}.
\end{equation}
%where $\vartheta=\frac{a}{b}$.
Without loss of generality, we first assume $m=6$ and $a=b=0.1$. Then $\vartheta=1$ and $\sqrt{\frac{1+\vartheta}{\vartheta}} \ \Big| \sqrt{a(1+\vartheta)^m}-1\Big|=2.1635<\sqrt{6}=2.4495$. Hence the condition \eqref{conn} is violated and no oscillation around $(a/b,a/b)=(1,1)$ is expected for the linearized system. To verify if this is the case for the nonlinear system \eqref{1.1}, we set the initial value $(u_0, v_0)$ as
\begin{equation}\label{ini1}
u_0(x)=v_0(x)=\frac{1}{1+e^{2(x-20)}}
\end{equation}
and perform the numerical simulations in an interval $[0, 200]$  with Dirichlet boundary conditions compatible with the initial value at the boundary. The numerical solution of \eqref{1.1} is shown in Fig.\ref{fig1} where we obverse that the solution will stabilize into monotone traveling waves although oscillates initially. This is also well consistent with  our analytical results about the existence of traveling wave solutions  given in Theorem \ref{main-thm1} when $\mathcal{K}(m,a)=0.4143<1$ if $m=6$ and $a=b=0.1$. Next we choose $m=4$ and $a=b=1$ such that $\sqrt{\frac{1+\vartheta}{\vartheta}} \ \Big| \sqrt{a(1+\vartheta)^m}-1\Big|=4.2426$ and hence \eqref{conn} holds. But numerically we still find that the system \eqref{1.1} will generate monotone traveling waves qualitatively similar to the patterns shown in Fig.\ref{fig1} (not shown here for brevity). This implies that the condition \eqref{conp} is not sufficient to induce non-monotone traveling waves oscillating around $(a/b, a/b)$.
%If we choose other forms of $\gamma(v)$ which decays without changing its convexity, we find that \eqref{1.1} generates similar monotone traveling wave solutions.

\begin{figure}[htbp]
\centering
\includegraphics[width=14cm]{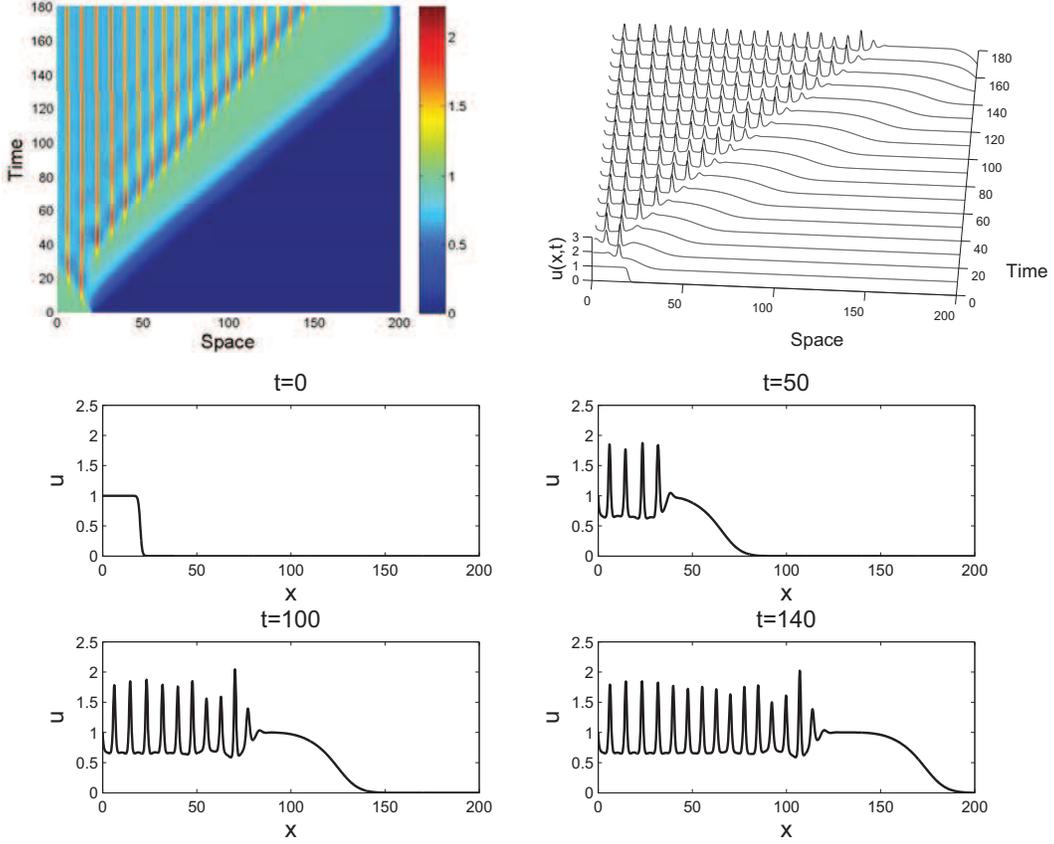}

\caption{Numerical simulations of wave propagation generated by the system \eqref{1.1} in $[0,200]$ with $\gamma(v)=1-\frac{v-1}{\sqrt{0.1+(v-1)^2}}, a=b=0.2, u_0=v_0=\frac{1}{1+e^{2(x-20)}}$.}
\label{fig2}
\end{figure}

Now an important question is whether the density-suppressed motility system \eqref{1.1} is capable of producing persistent oscillating traveling waves to interpret (at least qualitatively) the pattern observed in the experiment (see Fig.\ref{fig0}).  To explore this question numerically, we  consider the following sigmoid motility function
\begin{equation}\label{motility2}
\gamma(v)=1-\frac{v-1}{\sqrt{0.1+(v-1)^2}}
\end{equation}
which decays but changes the convexity at the point $v=1$, in contrast to the decreasing function \eqref{motility} whose convexity remains unchanged. We perform the numerical simulations for \eqref{1.1} with $a=b=0.2$ in an interval $[0,200]$ with the same initial value \eqref{ini1}. Remarkably we find non-monotone traveling wavefronts develop (see Fig.\ref{fig2}) and persist in time, where the wave oscillates at the trailing edge and propagates into the far field as time evolves. This is a prominent feature different from the patterns shown in Fig.\ref{fig1} generated from the motility function \eqref{motility}. If we choose some other forms of decreasing function $\gamma(v)$ that changes its convexity at $v=a/b=1$, we shall numerically find  similar non-monotone traveling wavefront patterns generated by \eqref{1.1}.

The above numerical simulations indicate, although not proved in this paper, that the density-suppressed motility system \eqref{1.1} can generate both monotone and non-monotone traveling wavefront solutions connecting $(a/b, a/b)$ to $(0,0)$. It numerically appears that the change of convexity of $\gamma(v)$  at $v=a/b$ is necessary to generate the non-monotone traveling wavefronts oscillating at the trailing edge around the equilibrium  $(a/b, a/b)$. The underlying mechanism remains mysterious and we will leave it as an open question for future study.

\begin{figure}[htbp]
\centering
\includegraphics[width=14cm]{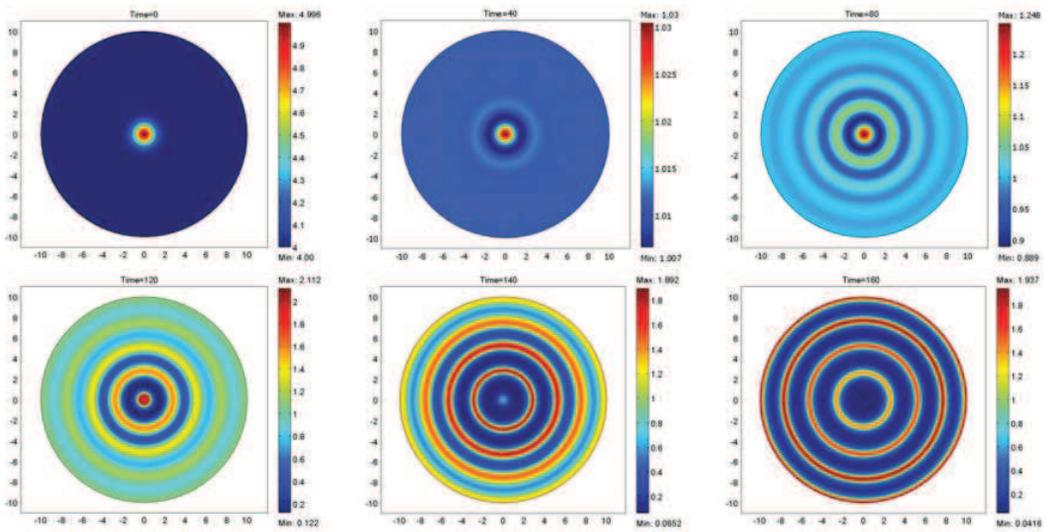}

\caption{Snapshot of numerical simulations of outward expanding ring patterns in a disk generated by the system \eqref{1.1} with $\gamma(v)=\frac{1}{(1+v)^6}, a=b=0.1, u_0=v_0=4+e^{-(x^2+y^2)}$.% $\gamma(v)=1-\frac{v-1}{\sqrt{0.1+(v-1)^2}}, a=b=0.1, u_0=v_0=5+e^{-(x^2+y^2)}$.
}
\label{fig3}
\end{figure}

Next we are devoted to exploring the patterns in a disk to mimic the apparatus used in the experiment of \cite{Liu-Science} where the experiment was conducted in petri dishes with bacteria initially inoculated at the center (see Fig.\ref{fig0}). In the numerical simulations, we set the domain as a disk with radius 10 and initially place the initial value $(u_0, v_0)=(4+e^{-(x^2+y^2)},4+e^{-(x^2+y^2)})$ in the center. We use the motility function given in \eqref{motility} with $m=6$ and set out Neumann boundary (i.e. zero-flux) conditions aligned with the experiment reality.  The snapshots of numerical patterns are recorded in Fig.\ref{fig3}, where we do observe the outward expanding ring patterns qualitatively analogous to the experiment patterns shown in Fig.\ref{fig0}. This validates the capability of model \eqref{1.1} reproducing the experimental patterns. However we should underline that it appears that the generation of oscillating patterns in two dimensions does not rely on the change of convexity of the motility function $\gamma(v)$ as shown in Fig.\ref{fig3}, which is very different from the situation in 1-D as shown in Fig.\ref{fig1} and Fig.\ref{fig2}. This imposes another interesting question  elucidating this subtle difference.
%As time evolves, the bacteria will run into boundary. Then patterns start to differentiate (at time $T=80$) due to the interaction of bacteria with the boundary  and evolves into other patterns consequently.

\subsection{Leading edge wave speeds} Following the spirit of classical method as in \cite{Mollison, Murray-book-2001}, we discuss the selection of the wave speed $c$ from the initial conditions given at infinity. Suppose that the initial value $(u_0, v_0)$ of the system \eqref{1.1} satisfies
\begin{equation}\label{ini}
\begin{cases}
u_0(x) \sim Ae^{-\lambda x}, \\
v_0(x) \sim Be^{-\lambda x},
\end{cases} \ \mathrm{as} \ x\ \to \infty
\end{equation}
with positive amplitudes $A$ and  $B$. Now we look for traveling wave solutions of \eqref{1.2} at the leading edge (i.e. $x \to \infty$) in the form of
 \begin{equation}\label{ansatz2}
\begin{cases}
u(x,t) \sim Ae^{-\lambda (x-ct)}, \\
v(x,t) \sim Be^{-\lambda (x-ct)}.
\end{cases}
\end{equation}
%where we have assumed that $u$ and $v$ have the same decay rates since $u \sim v$ as $x \to \infty$ from the second equation of \eqref{1.1}.
We substitute \eqref{ansatz2} into the first equation of \eqref{1.1} and get the dispersion relation between the wave speed $c$ and the initial decay rate $\lambda$:
\begin{equation}\label{speed1}
c=\gamma(0)\lambda+\frac{a}{\lambda}.
\end{equation}
Hence by the standard argument as in \cite{Murray-book-2001}, the asymptotic wave speed $c$ of traveling wave solutions to \eqref{1.1} satisfies
 \begin{equation}\label{speed2}
c=\begin{cases}
\gamma(0)\lambda +\frac{a}{\lambda}, & \text{if} \ \ 0<\lambda <\sqrt{a}, \\
2 \sqrt{\gamma(0)a}, & \text{if}\ \ \lambda \geq \sqrt{a}.
\end{cases}
\end{equation}
Next we plug \eqref{ansatz2} into the second equation of \eqref{1.1} and get the following relation on the amplitude of $u$ and $v$
\begin{equation}\label{amplitude}
A=[1+a+(\gamma(0)-1)\lambda^2]B.
\end{equation}
Therefore given the initial condition \eqref{ini}, the leading edge of traveling waves is fully determined by the ansatz \eqref{ansatz2} with wave speed \eqref{speed2} and amplitudes fulfilling \eqref{amplitude}.

As an example, we consider the motility function \eqref{motility} chosen in this paper, where $\gamma(0)=1$ and hence \eqref{speed1} gives
$$\lambda^2-c \lambda+a=0$$
which is exactly the same as the equation \eqref{1.23}. Furthermore \eqref{amplitude} gives $A=(1+a)B$ which well agrees with the result \eqref{am} in Theorem \ref{main-thm1}.

% the leading edge of evolving wave where $u$ is small and hence we can neglect $u^2$ and linearized equation of the first equation of \eqref{1.1} at $(u,v)=(0,0)$ is $$u_t=\gamma(0)u_{xx}+au.$$

\bigbreak
\noindent \textbf{Acknowledgment.}
 The research of Z.A. Wang was supported by the Hong Kong RGC GRF grant No. 15303019 (Project P0030816).

\end{document}